\documentclass[12pt]{article}
\usepackage{amsmath,amssymb,amsfonts,amsthm}

\topmargin- 0.5cm
 \small
{\par\vskip8pt minus3pt\rm}
\newcounter{item}[section]
\newcounter{kirshr}
\newcounter{kirsha}
\newcounter{kirshb}
\newenvironment{mysect}[1]{\vskip8pt\par\noindent\setcounter{item}{1}
\setcounter{equation}{0}
{\large\bf\arabic{section}.  #1 }\vskip8pt\nopagebreak\par\nopagebreak }
{\stepcounter{section}\upshape\par}
\newtheorem{theorem}{Theorem}[section]

\newtheorem{lemma}[theorem]{Lemma}
\newtheorem{corollary}[theorem]{Corollary}
\newtheorem{proposition}[theorem]{Proposition}

\newtheorem{remark}[theorem]{Remark}

\newtheorem{definition}[theorem]{Definition}
\newcommand\overcirc[1]{\raisebox{10pt}{\tiny$\circ$}{\kern-6.5pt}\mbox{$#1$}}
\newcommand\undersym[2]{\raisebox{-6pt}{\tiny$#2$}{\kern-5pt}\mbox{$#1$}}

\begin{document}
\title{\bf{Extended Absolute Parallelism Geometry}}

\author{Nabil. L. Youssef
 \ and  A. M.  Sid-Ahmed}

\date{}

\maketitle
\vspace{-1.13cm}

\begin{center}
{Department of Mathematics, Faculty of Science, \\Cairo
University, Giza, Egypt}\end{center} \vspace{-0.5cm}
\begin{center}
{nlyoussef2003@yahoo.fr, nyoussef@frcu.eun.eg\\ amrsidahmed@gmail.com, amrs@mailer.eun.eg}
\end{center}

\vspace{1cm} \maketitle
\smallskip

\noindent{\bf Abstract.} In this paper, we study Absolute
Parallelism (AP-) geometry on the tangent bundle $TM$ of a manifold
$M$. Accordingly, all geometric objects defined in this geometry are
not only functions of the positional argument $x$, but also depend
on the directional argument $y$. Moreover, many new geometric
objects, which have no counterpart in the classical AP-geometry,
emerge in this different framework. We refer to such a geometry as
an Extended Absolute Parallelism (EAP-) geometry.
\par
The building blocks of the EAP-geometry are a nonlinear connection
assumed given {\it a priori} and $2n$ linearly independent vector
fields (of special form) defined globally on $TM$ defining the
parallelization. Four different $d$-connections are used to explore
the properties of this geometry. Simple and compact formulae for the
curvature tensors and the W-tensors of the four defined
$d$-connections are obtained, expressed in terms of the torsion and
the contortion tensors of the EAP-space.
\par
Further conditions are imposed on the canonical $d$-connection
assuming that it is of Cartan type (resp. Berwald type). Important
consequences of these assumptions are investigated. Finally, a
special form of the canonical $d$-connection is studied under which
the classical AP-geometry is recovered naturally from the
EAP-geometry. Physical aspects of some of the geometric objects
investigated are pointed out and possible physical implications of
the EAP-space are discussed, including an outline of a generalized
field theory on the tangent bundle $TM$ of $M$. \footnote{ArXiv
number: 0805.1336}

\bigskip


\bigskip



\newpage

\begin{mysect}{Introduction}

The geometry of parallelizable  manifolds
or the Absolute Parallelism geometry (AP-geometry) (\cite{FI}, \cite {HP}, \cite{aa},
\cite{b}, \cite{AMR}) has many advantages in comparison
to Riemannian geometry. Unlike Riemannian geometry, which has ten
degrees of freedom (corresponding to the metric components for $n =
4$), AP-geometry has sixteen degrees of freedom (corresponding to
the number of components of the four vector fields defining the
parallelization). This makes AP-geometry a potential candidate for
describing physical phenomena other than gravity. Moreover, as
opposed to Riemannian geometry, which admits only one symmetric
linear connection, AP-geometry admits at least four natural
(built-in) linear connections, two of which are non-symmetric and
three of which have non-vanishing curvature tensors. Last, but not
least, associated with an AP-space there is a Riemannian structure
defined in a natural way. Thus, AP-geometry contains within its
geometric structure all the mathematical machinery of Riemannian
geometry. Accordingly, a comparison can be made between the results obtained in
the context of AP-geometry and general relativity, which is based on
Riemannian geometry.
\par
The geometry of the tangent bundle $(TM, \ \pi, \ M)$ of a smooth
manifold $M$ is very rich. It contains a lot of geometric objects of
theoretical interest and of a great importance in the construction
of various geometric models which have proved very useful in
different physical theories. Examples of such theories are the
general theory of relativity, particle physics, relativistic optics
and others.
\par
In this paper, we study AP-geometry in a context different from the classical one. Instead of
dealing with geometric objects defined on the manifold $M$, as in
the case of classical AP-space, we are dealing with geometric
objects defined on the tangent bundle $TM$ of $M$.
Accordingly, all geometric objects considered
are, in general, not only functions of the positional argument $x$, but also
depend on the directional argument $y$.
\par
The paper is organized in the following manner. In section 1,
following the introduction, we give a brief account of the basic
concepts and definitions that will be needed in the sequel. The
definitions of a $d$-connection, $d$-tensor field, torsion,
curvature, $hv$-metric and metric $d$-connection on $TM$ are
recalled. We end this section by the construction of a (unique)
metric $d$-connection on $TM$ which we refer to as the natural
metric $d$-connection. In section 2, we introduce the Extended
Absolute Parallelism (EAP-) geometry by assuming that $TM$ is
parallelizable \cite{B} and equipped with a nonlinear connection.
The canonical $d$-connection is then defined, expressed in terms of
the natural metric $d$-connection. In analogy to the classical
AP-geometry, two other $d$-connections are introduced: the dual and
the symmetric $d$-connections. We end this part with a comparison
between the classical AP-geometry and the EAP-geometry. In section
3, we carry out the task of computing the different curvature
tensors of the four defined $d$-connections. They are expressed, in
a relatively compact form, in terms of the torsion and the
contortion tensors of the EAP-space. All admissable contractions of
these curvature tensors are also obtained. In section 4, we
introduce and investigate the different $W$-tensors corresponding to
the different $d$-connections defined in the EAP-space, which are
again expressed in terms of the torsion and the contortion tensors.
In sections 5 and 6, we assume that the canonical $d$-connection is
of Cartan and Berwald type respectively. Some interesting results
are obtained, the most important of which is that, in the Cartan
type case, the given nonlinear connection is not independent of the
vector fields forming the parallelization, but can be expressed in
terms of their vertical counterparts. In section 7, we further
assume that the canonical $d$-connection is both of Cartan and
Berwald type. We show that, under this assumption, the classical
AP-geometry is recovered, in a natural way, from the EAP-geometry.
In section 8, we end this paper with some concluding remarks which
reveal possible physical applications of the EAP-space; among them
is an outline of a generalized field theory on the tangent bundle
$TM$ of $M$, based on Euler-Lagrange equations \cite{GLS} applied to
a suitable scalar Lagrangian.

\end{mysect}


\begin{mysect}{Fundamental Preliminaries}

In this section we give a brief account of the basic concepts and definitions that will be needed in the sequel.
Most of the material covered here may be found in \linebreak \cite{GLS}, \cite {V} with some slight modifications.

\bigskip

Let $M$ be a paracompact manifold of dimension $n$ of class $C^{\infty}$. Let $\pi:TM\to M$ be its tangent bundle.
If $(U, \ x^{\mu})$
is a local chart on $M$, then $(\pi^{-1}(U), (x^{\mu}, \ y^{a}))$ is the corresponding local
chart on $TM$.  The coordinate transformation on $TM$ is given by:
$$x^{\mu'} = x^{\mu'}(x^{\nu}), \ \ y^{a'} = p^{a'}_{a} y^{a},$$
$\mu = 1, ..., n; \ a = 1, ..., n$; \ $p^{a'}_{a} = \frac{\partial y^{a'}}{\partial y^{a}} =
\frac{\partial x^{a'}}{\partial x^{a}}$ and
det$(p^{a'}_{a})\neq 0$.
The paracompactness of $M$ ensures the existence of a nonlinear connection $N$ on $TM$
with coefficients $N^{a}_{\alpha}(x, y)$. The transformation formula for the coefficients
$N^{a}_{\alpha}$ is given by
\begin{equation}N^{a'}_{\alpha'} = p^{a'}_{a} p^{\alpha}_{\alpha'}N^{a}_{\alpha} + p^{a'}_{a}
p^{a}_{c'\alpha'}y^{c'},\end{equation}
where $p^{a}_{c'\alpha'} = \frac{\partial p^{a}_{c'}}{\partial x^{\alpha'}}$ .
The nonlinear connection leads to the direct sum decomposition
\begin{equation}\label{Sum}T_{u}(TM) = H_{u}(TM)\oplus V_{u}(TM), \ \ \forall u\in
TM\setminus \{0\}.\end{equation}
Here, $V_{u}(TM)$ is the vertical space at $u$ with
local basis $\dot \partial_a := \frac{\partial}{\partial y^{a}}$,
whereas $H_{u}(TM)$ is the horizontal space at $u$ associated with $N$
supplementary to the vertical space $V_{u}(TM)$. The canonical basis of $H_{u}(TM)$
is given by \begin{equation}\delta_{\mu}: = \partial_{\mu} - N^{a}_{\mu} \ \dot {\partial_a},\end{equation}
where
$\partial_{\mu}: = \frac{\partial}{\partial x^{\mu}}$.
Now, let $(dx^{\alpha}, \ \delta y^{a})$ be the basis of $T^{*}_{u}(TM)$ dual to the adapted
basis $(\delta_{\alpha}, \ \dot {\partial_a})$ of $T_{u}(TM)$. Then
\begin{equation}\delta y^{a} := dy^{a} + N^{a}_{\alpha}dx^{\alpha}\end{equation}
and
\begin{equation}dx^{\alpha}(\delta_{\beta}) = \delta^{\alpha}_{\beta}, \ \
dx^{\alpha}(\dot{\partial_{a}}) = 0; \ \ \ \delta y^{a}(\delta_{\beta}) = 0, \ \
\delta y^{a}(\dot \partial_{b}) = \delta^{a}_{b}.\end{equation}

\vspace{0.2cm}
Any vector field $X\in \mathfrak{X}(TM)$ is uniquely decomposed in the form $X = hX + vX$, where
$h$ and $v$ are respectively the horizontal and the vertical projectors associated with
the decomposition (\ref{Sum}). In the
adapted frame $(\delta_{\nu}, \ \dot{\partial_{a}})$, $hX = X^{\alpha}\delta_{\alpha}$ and $vX = X^{a}\dot {\partial_{a}}.$

\begin{definition} A nonlinear connection $N^{a}_{\mu}$ is said to be homogeneous if it is
positively homogeneous of degree $1$ in the directional argument $y$.
\end{definition}

\begin {definition} A $d$-connection on $TM$ is a linear connection on $TM$ which preserves by
parallelism the horizontal and vertical distribution:
if $Y$ is a horizontal (vertical) vector field, then $D_{X} Y$ is a horizontal (vertical) vector field, \  for all $X\in \mathfrak{X}(TM)$.
\end{definition}
Consequently, a $d$-connection
$D$ on $TM$ has only four coefficients. The coefficients
of a $d$-connection $D = (\Gamma^{\alpha}_{\mu\nu}, \ \Gamma^{a}_{b\nu}, \ C^{\alpha}_{\mu c}, \ C^{a}_{bc})$
are defined by
\begin{equation}D_{\delta \nu}\delta_{\mu} = :\Gamma^{\alpha}_{\mu\nu}\delta_{\alpha}, \ \ \
D_{\delta \nu}\dot{\partial_b} = :\Gamma^{a}_{b\nu}\dot {\partial_{a}}; \ \ \
D_{\dot{\partial_c}}\delta_{\mu} =: C^{\alpha}_{\mu c}\delta_{\alpha}, \ \ \
D_{\dot{\partial_c}}\dot {\partial_{b}} =: C^{a}_{bc}\dot{\partial_{a}}.\end{equation}

The transformation formulae of a $d$-connection are given by:
$$\Gamma^{\alpha'}_{\mu'\nu'} = p^{\alpha'}_{\alpha} p^{\mu}_{\mu'} p^{\nu}_{\nu'}
\Gamma^{\alpha}_{\mu\nu} + p^{\alpha'}_{\epsilon} p^{\epsilon}_{\mu'\nu'}, \ \ \ \Gamma^{a'}_{b'\mu'} = p^{a'}_{a} p^{b}_{b'} p^{\mu}_{\mu'}
\Gamma^{a}_{b\mu} + p^{a'}_{c}p^{c}_{b'\mu'};$$
$$C^{\alpha'}_{\mu' c'} = p^{\alpha'}_{\alpha} p^{\mu}_{\mu'} p^{c}_{c'}
C^{\alpha}_{\mu c}, \ \ \ C^{a'}_{b'c'} = p^{a'}_{a} p^{b}_{b'} p^{c}_{c'}
C^{a}_{bc}.$$

{\bf A comment on notation:} Both Greek indices $\{\alpha, \beta, \mu,...\}$ and
Latin indices $\{a, b, c, ...\}$, as previously mentioned,
take values from the same set $\{1,...,n\}$. It should be noted, however, that
{\bf Greek} indices are used to denote
{\bf horizontal} counterpart, whereas {\bf Latin} indices are used to denote {\bf vertical}
 counterpart. Einstein convention is applied on both types of indices.

\begin{definition} A $d$-tensor field $T$ on $TM$ of type $(p, r; \ q, s)$ is a tensor field on
$TM$ which can be locally expressed in the form
$$T = T^{{u_{1}}...{u_{p + r}}}_{{v_1}...{v_{q+s}}}\partial_{u_{1}}\otimes...
\otimes{\partial_{u_{p + r}}}\otimes dx^{v_1}\otimes... \otimes dx^{v_{q + s}},$$
where $u_{i}\in \{\alpha_{i}, \ a_i\}, \ v_{j}\in \{\beta_{j}, \ b_{j}\}$,
$$ \ \partial_{u_{i}}\in \{\delta_{\alpha_{i}}, \ \dot{\partial_{a_i}}\}, \
dx^{v_{j}}\in \{dx^{\beta_j}, \ \delta y^{b_j}\}, \ \ i = 1, ..., p + r; \ j = 1, ..., q + s,$$
so that the number of $\alpha_{i}$'s = $p$, the number of $a_{i}$'s = $r$, the number of
$\beta_{j}$'s = $q$ and the number of $b_{j}$'s = $s$.
\end{definition}

Let $T = T^{\alpha a}_{\beta b}\delta_{\alpha}\otimes{\dot{\partial_a}}\otimes dx^{\beta}
\otimes \delta y^{b}$ be a $d$-tensor field of type $(1, 1; 1, 1)$. Let $X\in \mathfrak{X}(TM)$ be such that
$X = hX + vX = X^{\mu}\delta_{\mu} + X^{c}\dot{\partial_c}$.
Then, by the properties of a $d$-connection, we have
$$D^{h}\!\!\,_{X} T: = D_{hX}T = (X^{\mu}T^{\alpha a}\!\,_{\beta b|\mu})
\delta_{\alpha}\otimes{\dot{\partial_a}}\otimes dx^{\beta}
\otimes \delta y^{b},$$
where
\begin{equation}\label{HCD}T^{\alpha a}\!\,_{\beta b|\mu}: = \delta_{\mu} T^{\alpha a}_{\beta b} +
T^{\epsilon a}_{\beta b}\Gamma^{\alpha}_{\epsilon \mu} +
T^{\alpha d}_{\beta b}\Gamma^{a}_{d\mu} -
T^{\alpha a}_{\epsilon b}\Gamma^{\epsilon}_{\beta \mu} -
T^{\alpha a}_{\beta d}\Gamma^{d}_{b\mu}.\end{equation}
Similarily,
$$D^{v}\!\!\,_{X} T: = D_{vX}T = (X^{c}T^{\alpha a}\!\,_{\beta b||c})
\delta_{\alpha}\otimes{\dot{\partial_a}}\otimes dx^{\beta}
\otimes \delta y^{b},$$
where
\begin{equation}\label{VCD}T^{\alpha a}\!\,_{\beta b||c}: = \dot{\partial_c} T^{\alpha a}_{\beta b} +
T^{\epsilon a}_{\beta b}C^{\alpha}_{\epsilon c} +
T^{\alpha d}_{\beta b}C^{a}_{d c} -
T^{\alpha a}_{\epsilon b}C^{\epsilon}_{\beta c} -
T^{\alpha a}_{\beta d}C^{d}_{bc}.\end{equation}

It is evident that (\ref{HCD}) and (\ref{VCD}) can be written for any
$d$-tensor field of arbitrary type.

\begin{definition} The two operators $D^{h}\!\!\,_{X}$ (denoted locally by $|$) and
$D^{v}\!\!\,_{X}$ (denoted locally by $||$) are called respectively the
horizontal ($h$-) and vertical ($v$-) covariant derivatives associated with the
$d$-connection $D$.
\end{definition}

\begin{definition} The torsion $\bf{T}$ of a d-connection $D$ on $TM$ is defined by
\begin{equation}\label{def}{\bf T}(X, Y): = D_{X} Y - D_{Y} X - [X, Y]; \ \ \forall X, Y\in \mathfrak{X}(TM).\end{equation}
\end{definition}
For getting the local expression for ${\bf T}$, we first recall that
$$[\delta_{\mu}, \ \delta_{\nu}] = R^{a}_{\mu\nu}\dot{\partial_a}; \ \ \
[\delta_{\mu}, \ \dot{\partial_b}] = (\dot{\partial_b} N^{a}_{\mu})\dot{\partial_a},$$
where
\begin{equation}R^{a}_{\mu\nu}: = \delta_{\nu}N^{a}_{\mu} - \delta_{\mu} N^{a}_{\nu}\end{equation}
is the
curvature of the nonlinear connection.

\bigskip

By a direct substitution in formula (\ref{def}), we obtain
\begin{proposition} In the adapted basis $(\delta_{\alpha}, \ \dot{\partial_{a}})$,
the torsion tensor ${\bf T}$ of a d-connection
$D = (\Gamma^{\alpha}_{\mu\nu}, \, \Gamma^{a}_{b \mu}, \, C^{\alpha}_{\mu c}, \, C^{a}_{bc})$
is charaterized by the following $d$-tensor fields with the local coefficients
$(\Lambda^{\alpha}_{\mu\nu}, \, R^{a}_{\mu\nu}, \, C^{\alpha}_{\mu c}, \, P^{a}_{\mu c}, \, T^{a}_{bc})$ defined by:
$$h{\bf T}(\delta_{\nu}, \ \delta_{\mu}) =: \Lambda^{\alpha}_{\mu\nu}\delta_{\alpha}, \ \ \ \ \
v{\bf T}(\delta_{\nu}, \ \delta_{\mu}) =:
R^{a}_{\mu\nu}\dot{\partial_a}$$
$$h{\bf T}(\dot{\partial_c}, \ \delta_{\mu}) =: C^{\alpha}_{\mu c}\delta_{\alpha}, \ \ \ \ \
v{\bf T}(\dot{\partial_c}, \ \delta_{\mu}) =: P^{a}_{\mu c}\dot{\partial_a}, \ \ \ \ \ \ v{\bf T}(\dot{\partial_c}, \
\dot{\partial_b}) =: T^{a}_{bc}\dot{\partial_a},$$
where
\begin{equation}\Lambda^{\alpha}_{\mu\nu} := \Gamma^{\alpha}_{\mu\nu} - \Gamma^{\alpha}_{\nu\mu},
\ \ \ \ P^{a}_{\mu c}: = \dot{\partial_c} N^{a}_{\mu} - \Gamma^{a}_{c \mu}, \ \ \ \ T^{a}_{bc} := C^{a}_{bc} - C^{a}_{cb}.
\end{equation}
\end{proposition}

Throughout the paper we shall use the notation
${\bf T} = (\Lambda^{\alpha}_{\mu\nu}, \, R^{a}_{\mu\nu}, \, C^{a}_{\mu c}, \, P^{a}_{\mu c}, \, T^{a}_{bc})$.

\begin{corollary} The torsion tensor ${\bf T} = (\Lambda^{\alpha}_{\mu\nu}, \, R^{a}_{\mu\nu}, \, C^{a}_{\mu c}, \, P^{a}_{\mu c}, \, T^{a}_{bc})$
of a $d$-connection $D$ vanishes if
$$\Gamma^{\alpha}_{\mu\nu} = \Gamma^{\alpha}_{\nu\mu}, \ \ \ R^{a}_{\mu\nu} = C^{\alpha}_{\mu c} = 0, \ \
\ \dot{\partial_c} N^{a}_{\mu} =  \Gamma^{a}_{c \mu}, \ \ \ C^{a}_{bc} = C^{a}_{cb}.$$
\end{corollary}

\begin{definition}The curvature tensor ${\bf R}$ of a $d$-connection $D$ is given by
$${\bf R}(X, Y)Z: = D_{X} D_{Y} Z - D_{Y}D_{X} Z - D_{[X, \ Y]}Z; \ \ \forall X, Y, Z\in \mathfrak{X}(TM).$$
\end{definition}

By definition of a $d$-connection, it follows that ${\bf R}(X, Y)Z$ is determined by eight $d$-tensor fields,
six of which are independent
due to the fact that
${\bf R}(X, Y) = - {\bf R}(Y, X)$. We set
$$ \ {\bf R}(\delta_{\mu}, \ \delta_{\nu})\delta_{\beta} =: R^{\alpha}_{\beta\mu\nu}\delta_{\alpha};
\ \ \ \ \ {\bf R}(\delta_{\mu}, \ \delta_{\nu})\dot{\partial_{b}} =: R^{a}_{b\mu\nu}\dot{\partial_{a}},$$
$${\bf R}(\dot{\partial_c}, \ \delta_{\nu})\delta_{\beta} =: P^{\alpha}_{\beta\nu c}\delta_{\alpha};
\ \ \ \ \ \ {\bf R}(\dot{\partial_c}, \ \delta_{\nu})\dot{\partial_{b}} =: P^{a}_{b\nu c}\dot{\partial_{a}},$$
$${\bf R}(\dot{\partial_b}, \ \dot{\partial_c})\delta_{\beta} =: S^{\alpha}_{\beta bc}\delta_{\alpha};
\ \ \ \ \ \ {\bf R}(\dot{\partial_c}, \ \dot{\partial_d})\dot{\partial_{b}} =: S^{a}_{bcd}\dot{\partial_{a}}.$$

Throughout the paper we shall use the notation ${\bf {R}} =
({R}^{\alpha}_{\beta\mu\nu}, \!\, {R}^{a}_{b \mu\nu}, \!\, {P}^{\alpha}_{\beta\nu c}, \!\, {P}^{a}_{b \nu c},
\!\, {S}^{\alpha}_{\beta bc}, \!\, {S}^{a}_{bcd})$.
\begin{theorem} \label{CUR} The curvature ${\bf R}$ of a $d$-connection
$D = (\Gamma^{\alpha}_{\mu\nu}, \ \Gamma^{a}_{b\mu}, \ C^{\alpha}_{\mu c}, \ C^{a}_{bc})$
is charaterized by the $d$-tensor fields with local coefficients:
\begin{description}
\item [(a)] $R^{\alpha}_{\beta\mu\nu} = \delta_{\mu}\Gamma^{\alpha}_{\beta\nu} -
\delta_{\nu}\Gamma^{\alpha}_{\beta\mu} + \Gamma^{\epsilon}_{\beta\nu}
\Gamma^{\alpha}_{\epsilon\mu} -  \Gamma^{\epsilon}_{\beta\mu}
\Gamma^{\alpha}_{\epsilon\nu} + C^{\alpha}_{\beta d}R^{d}_{\nu\mu},$
\item [(b)] $R^{a}_{b\mu\nu} = \delta_{\mu}\Gamma^{a}_{b\nu} -
\delta_{\nu}\Gamma^{a}_{b\mu} + \Gamma^{c}_{b\nu}
\Gamma^{a}_{c\mu} -  \Gamma^{c}_{b\mu}
\Gamma^{a}_{c\nu} + C^{a}_{bd}R^{d}_{\nu\mu},$
\item [(c)] $P^{\alpha}_{\beta\nu c} = \dot{\partial_c}\Gamma^{\alpha}_{\beta\nu} -
C^{\alpha}_{\beta c|\nu} + C^{\alpha}_{\beta d}P^{d}_{\nu c},$
\item [(d)] $P^{a}_{b\nu c} = \dot{\partial_c}\Gamma^{a}_{b\nu} -
C^{a}_{bc|\nu} + C^{a}_{bd}P^{d}_{\nu c},$
\item [(e)] $S^{\alpha}_{\beta bc} = \dot{\partial_b}C^{\alpha}_{\beta c} -
\dot{\partial_c}C^{\alpha}_{\beta b} + C^{\epsilon}_{\beta c}C^{\alpha}_{\epsilon b} -
C^{\epsilon}_{\beta b}C^{\alpha}_{\epsilon c},$
\item [(f)] $S^{a}_{bcd} = \dot{\partial_c}C^{a}_{bd} -
\dot{\partial_d}C^{a}_{bc} + C^{e}_{bd}C^{a}_{ec} -
C^{e}_{bc}C^{a}_{ed}.$
\end{description}
\end{theorem}

\begin{corollary}\label{vanish} The curvature tensor ${\bf {R}} =
({R}^{\alpha}_{\beta\mu\nu}, \!\, {R}^{a}_{b \mu\nu}, \!\, {P}^{\alpha}_{\beta\nu c}, \!\, {P}^{a}_{b \nu c},
\!\, {S}^{\alpha}_{\beta bc}, \!\, {S}^{a}_{bcd})$ of a $d$-connection $D$ vanishes iff
$$R^{\alpha}_{\beta\mu\nu} = R^{a}_{b\mu\nu} =
P^{\alpha}_{\beta\nu c} = P^{a}_{b\nu c} = S^{\alpha}_{\beta bc} = S^{a}_{bcd} = 0.$$
\end{corollary}

\begin{definition} An $hv$-metric on $TM$ is a covariant $d$-tensor field
${\bf G} := h{\bf G} + v{\bf G}$ on $TM$, where $h{\bf G}: =
g_{\alpha\beta}\,dx^{\alpha}\otimes dx^{\beta}$, \ $v{\bf G} :=
g_{ab}\,\delta y^{a}\otimes \delta y^{b}$ such that:
\begin{equation}g_{\alpha\beta} = g_{\beta\alpha}, \ \ det (g_{\alpha\beta})\neq 0; \ \
\ g_{ab} = g_{ba}, \ \ det (g_{ab})\neq 0.\end{equation}
\end{definition}
The inverses of $(g_{\alpha\beta})$ and $(g_{ab})$, denoted by $(g^{\alpha\beta})$ and $(g^{ab})$ repectively, are
given by
\begin{equation}g_{\alpha\epsilon}\,g^{\epsilon\beta} = \delta^{\alpha}_{\beta}, \ \ \ g_{ae}\,g^{eb} = \delta^{a}_{b}.\end{equation}

\begin{definition} A $d$-connection $D$ on $TM$ is said to be metric or compatible with
the metric ${\bf G}$ if $D_{X}{\bf G} = 0, \ \forall X\in \mathfrak{X}(TM)$.
\end{definition}

In the adapted frame $(\delta_{\alpha}, \ \dot{\partial_a})$, the above condition can
be expressed locally in the form:
\begin{equation}\label{locally}g_{\alpha\beta|\mu} = g_{\alpha\beta||c} = g_{ab|\mu} = g_{ab||c} = 0.\end{equation}

We have the following Theorem \cite{GLS}:
\begin{theorem}\label{metric} There exists a unique metrical $d$-connection \ $\overcirc{D} =
( \ \overcirc{\Gamma^{\alpha}_{\mu\nu}}, \ \, \overcirc{\Gamma^{a}_{b\nu}}, \
\, \overcirc{C^{\alpha}_{\mu c}}, \ \, \overcirc{C^{a}_{bc}})$ on $TM$ with the properties that
\begin{description}
\item [(a)] $\overcirc{\Lambda}^{\alpha}_{\mu\nu} = \ \overcirc{\Gamma}^{\alpha}_{\mu\nu} - \
\overcirc{\Gamma}^{\alpha}_{\nu\mu} = 0, \ \, \ \ \ \ \overcirc{T}^{a}_{bc} = \ \overcirc{C}^{a}_{bc} - \ \overcirc{C}^{a}_{cb} = 0.$
\item [(b)] $\overcirc{\Gamma}^{a}_{b\nu} := \dot{\partial_b} N^{a}_{\nu} +
\frac{1}{2} \, g^{ac}(\delta_{\nu} g_{bc} - g_{dc} \ \dot{\partial_b}N^{d}_{\nu} -
g_{bd} \ \dot{\partial_c}N^{d}_{\nu}), \ \ \ \ \ \ \, \overcirc{C}^{\alpha}_{\mu c} := \frac{1}{2} \ g^{\alpha\epsilon}\dot{\partial_c}g_{\mu\epsilon}.$
\end{description}
In this case, the coefficients \ $\overcirc{\Gamma}^{\alpha}_{\mu\nu}$ and \ $\overcirc{C}^{a}_{bc}$ are necessarily
of the form
$$\overcirc{\Gamma}^{\alpha}_{\mu\nu}: = \frac{1}{2} \ g^{\alpha\epsilon}(\delta_{\mu}
g_{\epsilon\nu} + \delta_{\nu}g_{\epsilon\mu} - \delta_{\epsilon}g_{\mu\nu}), \ \ \ \ \
\overcirc{C^{a}_{bc}}: = \frac{1}{2} \ g^{ad}(\dot{\partial_{b}}g_{dc} +
\dot{\partial_{c}}g_{db} - \dot{\partial_{d}}g_{bc}).$$
\end{theorem}

\begin{definition} The $d$-connection \
$\overcirc{D} = ( \ \overcirc{\Gamma}^{\alpha}_{\mu\nu}, \ \, \overcirc{\Gamma}^{a}_{b\nu},
\ \, \overcirc{C}^{\alpha}_{\mu c}, \ \, \overcirc{C}^{a}_{bc})$ defined in Theorem \ref{metric}
will be referred to as the {\bf natural metric} $d$-connection. The $h$- and $v$-covariant derivatives with respect to the natural
metric $d$-connection \ $\overcirc{D}$ will be denoted by ${o\atop|}$ and ${o\atop||}$ respectively.

\end{definition}
\end{mysect}


\begin{mysect}{Extended Absolute Parallelism Geometry (EAP-geometry)}

In this section, we study AP-geometry in a context different from the classical one. Instead of dealing with geometric objects
defined on the manifold $M$, we will be dealing with geometric objects defined on the tangent bundle $TM$ of $M$.
Many new geometric objects, which have no
counterpart in the classical AP-geometry, emerge in this different framework. Moreover, the
basic geometric objects of the new geometry acquire a richer structure compared to the
corresponding basic geometric objects of the classical AP-geometry (See Table 2).

\bigskip

As in the previous section, $M$ is assumed to be a smooth paracompact manifold of dimension $n$. This insures the existence of a
nonlinear connection on
$TM$ so that the decomposition (\ref{Sum}) induced by the nonlinear connection holds.

\bigskip

We assume that \ $\undersym{\lambda}{i}$, $i = 1, ..., n$, are $n$ vector fields defined {\bf globally}
on $TM$. In the adapted basis $(\delta_{\alpha}, \ \dot{\partial_{a}})$, we have \
$\undersym{\lambda}{i} = \ h \, \undersym{\lambda}{i} + \ v \, \undersym{\lambda}{i} =
 \ \undersym{\lambda}{i}^{\alpha}{\delta_{\alpha}} +
\ \undersym{\lambda}{i}^{a}\dot{\partial_a}.$ We further assume that the $n$ horizontal vector fields \
$h \, \undersym{\lambda}{i}$ and
the $n$ vertical vector fields \ $v \, \undersym{\lambda}{i}$ are {\bf linearly
independent}. This implies, in particular, that the $n$ vector fields \ $\undersym{\lambda}{i}$, themselves, are linearly independent. Moreover,
we have
\begin{equation}\label{inverse}\undersym{\lambda}{i}^{\alpha} \ \undersym{\lambda}{i}_{\beta} = \delta^{\alpha}_{\beta}, \
\ \ \ \undersym{\lambda}{i}^{\alpha} \ \undersym{\lambda}{j}_{\alpha} = \delta_{ij}; \ \ \ \
\ \undersym{\lambda}{i}^{a} \ \undersym{\lambda}{i}_{b} = \delta^{a}_{b}, \
\ \ \ \undersym{\lambda}{i}^{a} \ \undersym{\lambda}{j}_{a} = \delta_{ij},\end{equation}
where $( \ \undersym{\lambda}{i}_{\alpha})$ and $( \ \undersym{\lambda}{i}_{a})$ denote the
inverse matrices of $( \ \undersym{\lambda}{i}^{\alpha})$ and $( \ \undersym{\lambda}{i}^{a})$
respectively.

\bigskip

We refer to the above space, which we denote by $(TM, \ \lambda)$, as an
{\bf Extended Absolute Parellelism} (EAP-) geometry which is characterized by the
{\bf existence of $2n$ linearly independent vector fields defined globally on $TM$}.

\bigskip

The Latin indices
$\{i, j\}$ will be used for numbering the $n$ vector fields (mesh indices).
Einstein convention is applied on the mesh indices (which will
always be written in  lower position) as well as the component indices.
In the sequel, to simplify notations,
we will use the symbol $\lambda$ without
the subscript $i$ to denote any one of the vector fields \
$\undersym{\lambda}{i} \ (i = 1, ..., n)$. The index $i$
will appear {\it only when summation is performed.}

\bigskip

Let us define
\begin{equation}\label{apm}g_{\alpha\beta}: = \ \undersym{\lambda}{i}_{\alpha} \ \undersym{\lambda}{i}_{\beta}, \ \
\ \ g_{ab}: = \ \undersym{\lambda}{i}_{a} \ \undersym{\lambda}{i}_{b}.\end{equation} Then, clearly,
$${\bf G} =  g_{\alpha\beta} \,dx^{\alpha}\otimes dx^{\beta} + g_{ab}\,\delta y^{a}\otimes \delta y^{b}$$
is an $hv$-metric on $TM$. Moreover, in view of (\ref{inverse}), the
inverse of the matrices $(g_{\alpha\beta})$ and $(g_{ab})$ are given
by $(g^{\alpha\beta})$ and $(g^{ab})$ respectively, where
\begin{equation}g^{\alpha\beta} = \ \undersym{\lambda}{i}^{\alpha} \
\undersym{\lambda}{i}^{\beta}, \ \ \ \ g^{ab} = \ \undersym{\lambda}{i}^{a} \
\undersym{\lambda}{i}^{b}.\end{equation}


Now, let \ $\overcirc{D} = ( \ \overcirc{\Gamma}^{\alpha}_{\mu\nu}, \ \, \overcirc{\Gamma}^{a}_{b\nu},
\ \, \overcirc{C}^{\alpha}_{\mu c}, \ \, \overcirc{C}^{a}_{bc})$ be the natural metric $d$-connection defined by Theorem
\ref{metric}, where $g_{\mu\nu}$ and $g_{ab}$ are the metric tensors given by (\ref{apm}).
\begin{theorem}\label{apc} There exists a unique $d$-connection $D = (\Gamma^{\alpha}_{\mu\nu}, \
\Gamma^{a}_{b\nu}, \ C^{\alpha}_{\mu c}, \ C^{a}_{bc})$ such that
\begin{equation}\label{NY} {\lambda}^{\alpha}\!\, _{|\mu} = {\lambda}^{\alpha}\!\, _{||c} =
{\lambda}^{a}\!\, _{|\mu} = {\lambda}^{a}\!\, _{||c} = 0,\end{equation}
where $|$ and $||$ are the h- and v-covariant derivatives with respect to $D$.
Consequently $D$ is a metric $d$-connection. It is given by
\begin{equation}\label{cm}\Gamma^{\alpha}_{\mu\nu}: = \ \overcirc{\Gamma}^{\alpha}_{\mu\nu} +
\ \undersym{\lambda}{i}^{\alpha} \ \undersym{\lambda}{i}_{\mu{o\atop|}\nu}, \ \ \
\ \Gamma^{a}_{b\nu} := \ \overcirc{\Gamma}^{a}_{b\nu} + \ \undersym{\lambda}{i}^{a} \
\undersym{\lambda}{i}_{b{o\atop|}\nu};\end{equation}
\begin{equation}\label{ccm}C^{\alpha}_{\mu c} := \ \overcirc{C}^{\alpha}_{\mu c} +
\ \undersym{\lambda}{i}^{\alpha} \ \undersym{\lambda}{i}_{\mu{o\atop||}c}, \ \ \
\ C^{a}_{bc}: = \ \overcirc{C}^{a}_{bc} + \ \undersym{\lambda}{i}^{a} \ \undersym{\lambda}{i}_{b{o\atop||}c}.
\end{equation}

Relation (\ref{NY}) will be called the AP-condition (as in the classical AP-geometry). \end{theorem}

\begin{proof} First, it is clear that $D$ is a $d$-connection on $TM$.
We next prove that ${\lambda}^{\alpha}\!\, _{|\nu} = 0$.
We have
\begin{eqnarray*} {\lambda}^{\alpha}\!\,_{|\nu} &=& \delta_{\nu} \
{\lambda}^{\alpha} + {\lambda}^{\mu}
\Gamma^{\alpha}_{\mu\nu} = \delta_{\nu}{\lambda}^{\alpha} + \
\undersym{\lambda}{i}^{\mu}(\ \overcirc{\Gamma^{\alpha}_{\mu\nu}} +
\undersym{\lambda}{j}^{\alpha} \ \undersym{\lambda}{j}_{\mu{o\atop|}\nu})\\
&=& (\delta_{\nu}{\lambda}^{\alpha} + \ \overcirc{\Gamma^{\alpha}_{\mu\nu}}
{\lambda}^{\mu})
- ( \ \undersym{\lambda}{i}^{\mu} \ \undersym{\lambda}{j}_{\mu}) \ \undersym{\lambda}{j}
^{\alpha}\!\,_{{o\atop|}\nu} = {\lambda}^{\alpha}\!\,_{{o\atop|}\nu} - \
{\lambda}^{\alpha}\!\,_{{o\atop|}\nu} = 0.\\
\\[- 1.2 cm]\end{eqnarray*}
The rest is proved in a similar manner.
\end{proof}
\begin{definition} The $d$-connection $D = (\Gamma^{\alpha}_{\mu\nu}, \ \Gamma^{a}_{b\nu},
\ C^{\alpha}_{\mu c}, \ C^{a}_{bc})$ defined in Theorem \ref{apc} will be referred to as the {\bf canonical} $d$-connection
of the EAP-space.
\end{definition}

In analogy to the classical AP-space, the torsion tensor of the canonical $d$-connection will be called the torsion tensor of the EAP-space.
\begin{theorem}\label{canonical} The canonical $d$-connection $D$ can be expressed explicitely in terms of
the ${\lambda}$'s only in the form:
\begin{equation}\Gamma^{\alpha}_{\mu\nu} =
\ \undersym{\lambda}{i}^{\alpha}(\delta_{\nu} \ \undersym{\lambda}{i}_{\mu}), \ \ \ \ \ \
\Gamma^{a}_{b\nu} = \ \undersym{\lambda}{i}^{a}(\delta_{\nu} \ \undersym{\lambda}{i}_{b});\end{equation}
\begin{equation}C^{\alpha}_{\mu c} = \ \undersym{\lambda}{i}^{\alpha} (\dot{\partial_c} \
\undersym{\lambda}{i}_{\mu}), \ \ \ \ \ \ C^{a}_{bc} = \ \undersym{\lambda}{i}^{a} (\dot{\partial_c} \
\undersym{\lambda}{i}_{b}).\end{equation}
\end{theorem}

\begin{proof} Since ${\lambda}^{\alpha}\!\, _{|\nu} = 0$, it follows that
$\delta_{\nu}{\lambda}^{\alpha} = - {\lambda}^{\epsilon}\Gamma^{\alpha}_{\epsilon\nu}$.
Multiplying by \!\, $\undersym{\lambda}{i}_{\mu}$, we get \
$\undersym{\lambda}{i}_{\mu} (\delta_{\nu} \ \undersym{\lambda}{i}^{\alpha}) = - \ \Gamma^{\alpha}_{\mu\nu}$
so that, by (\ref{inverse}), $\Gamma^{\alpha}_{\mu\nu} = \ \undersym{\lambda}{i}^{\alpha}(\delta_{\nu} \ \undersym{\lambda}{i}_{\mu})$.
The other formulae are derived in a similar manner.
\end{proof}

By Theorem \ref{apc} and Theorem \ref{canonical}, we have
\begin{corollary} The natural metric $d$-connection \ $\overcirc{D}$ can be expressed
explicitely in terms
of the $\lambda$'s only in the form
\begin{equation}\overcirc{\Gamma}^{\alpha}_{\mu\nu} = \ \undersym{\lambda}{i}^{\alpha}
(\delta_{\nu} \ \undersym{\lambda}{i}_{\mu} -\  \undersym{\lambda}{i}_{\mu{o\atop{|}}\nu}), \ \
\ \ \ \ \overcirc{\Gamma}^{a}_{b\nu} = \ \undersym{\lambda}{i}^{a}
(\delta_{\nu} \ \undersym{\lambda}{i}_{b} -\  \undersym{\lambda}{i}_{b{o\atop{|}}\nu});\end{equation}
\begin{equation}\overcirc{C}^{\alpha}_{\mu c} = \ \undersym{\lambda}{i}^{\alpha}
(\dot{\partial_c} \ \undersym{\lambda}{i}_{\mu} -\  \undersym{\lambda}{i}_{\mu{o\atop{||}}c}),
\ \ \ \ \ \ \overcirc{C}^{a}_{bc} = \ \undersym{\lambda}{i}^{a}
(\dot{\partial_c} \ \undersym{\lambda}{i}_{b} -\  \undersym{\lambda}{i}_{b{o\atop{||}}c}).\end{equation}
\end{corollary}

\begin{definition} The contortion tensor of an EAP-space is defined by
$${\bf C}(X, \ Y): = D_{Y}X - \ \overcirc{D}_{Y}X; \ \ \forall X, Y\in \mathfrak{X}(TM),$$
where $D$ is the cannonical $d$-connection and \ $\overcirc{D}$ is the natural metric
$d$-connection.
\end{definition}

In the adapted basis $(\delta_{\mu}, \ \dot{\partial_a})$, the contortion tensor
is characterized by the following $d$-tensor fields:
$${\bf C}(\delta_{\mu}, \ \delta_{\nu}) = :\gamma^{\alpha}_{\mu\nu}\delta_{\alpha}, \ \ \
{\bf C}(\delta_{\mu}, \ \dot{\partial_c}) =: \gamma^{\alpha}_{\mu c}\delta_{\alpha};
\ \ \ \ \ {\bf C}(\dot{\partial_b}, \ \delta_{\mu}) = :\gamma^{a}_{b \mu}\dot{\partial_a},
\ \ \ \ {\bf C}(\dot{\partial_{b}}, \ \dot{\partial_c}) =: \gamma^{a}_{b c}\dot{\partial_a};$$
\begin{equation}\label{cnt}\gamma^{\alpha}_{\mu\nu} := \Gamma^{\alpha}_{\mu\nu} -
\ \overcirc{\Gamma}^{\alpha}_{\mu\nu}, \ \ \ \gamma^{a}_{b\mu}: =
\Gamma^{a}_{b\mu} -
\ \overcirc{\Gamma}^{a}_{b\mu}; \ \ \ \ \gamma^{\alpha}_{\mu c}: = C^{\alpha}_{\mu c} -
\ \overcirc{C}^{\alpha}_{\mu c}, \ \ \ \gamma^{a}_{bc}: = C^{a}_{bc} - \ \overcirc{C}^{a}_{bc}.\end{equation}

Throughout the paper we shall use the notation ${\bf {C}} = (\gamma^{\alpha}_{\mu\nu}, \
\gamma^{a}_{b \mu}, \ \gamma^{\alpha}_{\mu c}, \ \gamma^{a}_{bc})$.

\bigskip

By definition of the canonical $d$-connection and (\ref{cnt}), the contortion tensor can be expressed explicitely in terms of
the ${\lambda}$'s only in the form:
\begin{equation}\label{contortion}\gamma^{\alpha}_{\mu\nu} = \ \undersym{\lambda}{i}^{\alpha} \
\undersym{\lambda}{i}_{\mu{o\atop{|}}\nu}, \ \
\gamma^{a}_{b \mu} = \ \undersym{\lambda}{i}^{a} \
\undersym{\lambda}{i}_{b {o\atop{|}}\mu}, \ \ \gamma^{\alpha}_{\mu c} =
\ \undersym{\lambda}{i}^{\alpha} \
\undersym{\lambda}{i}_{\mu{o\atop{||}}c}, \ \ \ \gamma^{a}_{b c} =
\ \undersym{\lambda}{i}^{a} \ \undersym{\lambda}{i}_{b {o\atop{||}}c}.\end{equation}

\begin{proposition}\label{skew} Let $\gamma_{\alpha\mu\nu}: =
g_{\alpha\epsilon}\gamma^{\epsilon}_{\mu\nu}, \
\gamma_{ab\mu}: = g_{ac}\gamma^{c}_{b \mu}, \ \gamma_{\alpha\mu c} :=
g_{\alpha\epsilon}\gamma^{\epsilon}_{\mu c}$, \ $\gamma_{abc} := g_{ad}\gamma^{d}_{bc}$.
Then each of the above defined $d$-tensor fields is skew-symmetric in the
first pair of indices.
Consequently,
$\gamma^{\alpha}_{\alpha \nu} = \gamma^{a}_{a\mu} = \gamma^{\alpha}_{\alpha c} =
\gamma^{a}_{a c} = 0.$
\end{proposition}

\begin{proof} We have
$$\gamma_{\alpha\mu\nu} + \gamma_{\mu\alpha\nu} =
\ \undersym{\lambda}{i}_{\alpha} \ \undersym{\lambda}{i}_{\mu{o\atop{|}}\nu} +
\ \undersym{\lambda}{i}_{\mu} \ \undersym{\lambda}{i}_{\alpha{o\atop{|}}\nu} =
( \ \undersym{\lambda}{i}_{\alpha} \ \undersym{\lambda}{i}_{\mu})_{{o\atop{|}}\nu} =
g_{\alpha\mu{o\atop{|}}\nu} = 0.$$
The rest is proved analogously.
\end{proof}

A simple calculation gives

\begin{proposition} Let ${\bf T} = (\Lambda^{\alpha}_{\mu\nu}, \
R^{a}_{\mu\nu}, \ C^{\alpha}_{\mu c}, \ P^{a}_{\mu b}, \ T^{a}_{bc})$ and ${\bf {C}} = (\gamma^{\alpha}_{\mu\nu}, \
\gamma^{a}_{b \mu}, \ \gamma^{\alpha}_{\mu c}, \ \gamma^{a}_{bc})$ be the
torsion and the contortion tensors of the EAP-space
respectively. Then the following relations hold:
\begin{equation}\label{torsion}\Lambda^{\alpha}_{\mu\nu} =  \gamma^{\alpha}_{\mu\nu} - \gamma^{\alpha}_{\nu\mu}, \ \
\ P^{a}_{\mu b} = - \gamma^{a}_{b \mu} + \ \overcirc{P}^{a}_{\mu b}, \ \ \
C^{\alpha}_{\mu c} = \gamma^{\alpha}_{\mu c} + \ \overcirc{C}^{\alpha}_{\mu c},
\ \ \ T^{a}_{bc} = \gamma^{a}_{bc} - \gamma^{a}_{cb}.\end{equation}
Consequently,
\begin{equation}\label{conT}\Lambda^{\alpha}_{\mu\alpha} = \gamma^{\alpha}_{\mu\alpha} = :C_{\mu},  \ \
T^{a}_{ba} = \gamma^{a}_{ba} = :C_{b}.\end{equation}
\end{proposition}

\begin{remark}\em{\label{TC}It can be shown, in analogy to the classical AP-space \cite{H},
that
\begin{equation}\label{tc}\gamma_{\alpha\mu\nu} = \frac{1}{2}(\Lambda_{\alpha\mu\nu} + \Lambda_{\nu\mu\alpha} +
\Lambda_{\mu\nu\alpha}), \ \ \ \gamma_{abc} = \frac{1}{2}(T_{abc} +
T_{cba} + T_{bca});\end{equation} where $\Lambda_{\alpha\mu\nu}: =
g_{\alpha\epsilon}\,\Lambda^{\epsilon}_{\mu\nu}$ and $T_{abc}: =
g_{ad}\,T^{d}_{bc}.$}
\end{remark}
By (\ref{torsion}) and (\ref{tc}), $\Lambda^{\alpha}_{\mu\nu}$ (resp. $T^{a}_{bc}$)
vanishes iff $\gamma^{\alpha}_{\mu\nu}$ (resp. $\gamma^{a}_{bc}$) vanishes.

\begin{definition} Let $D = (\Gamma^{\alpha}_{\mu\nu}, \ \Gamma^{a}_{b\mu}, \
C^{\alpha}_{\mu c}, \ C^{a}_{bc})$ be the canonical $d$-connection.

\begin{description}
\item [(a)] The {\bf dual} $d$-connection
$\widetilde{D} = (\widetilde{\Gamma}^{\alpha}_{\mu\nu}, \ \widetilde{\Gamma}^{a}_{b\mu}, \
\widetilde{C}^{\alpha}_{\mu c}, \ \widetilde{C}^{a}_{bc})$ is defined by
\begin{equation}\widetilde{\Gamma}^{\alpha}_{\mu\nu}: = \Gamma^{\alpha}_{\nu\mu}, \ \ \ \
\widetilde{\Gamma}^{a}_{b\mu}: = \Gamma^{a}_{b\mu}; \ \ \ \ \ \ \widetilde{C}^{\alpha}_{\mu c} :=
C^{\alpha}_{\mu c}, \ \ \ \ \widetilde{C}^{a}_{bc}: = C^{a}_{cb}.\end{equation}

\item [(b)] The {\bf symmertic} $d$-connection $\widehat{D} = (\widehat{\Gamma}^{\alpha}_{\mu\nu}, \
\widehat{\Gamma}^{a}_{b\mu}, \
\widehat{C}^{\alpha}_{\mu c}, \ \widehat{C}^{a}_{bc})$ is defined by
\begin{equation}\widehat{\Gamma}^{\alpha}_{\mu\nu}: = \frac{1}{2}(\Gamma^{\alpha}_{\mu\nu} +
\Gamma^{\alpha}_{\nu\mu}), \ \ \ \widehat{\Gamma}^{a}_{b\mu}: = \Gamma^{a}_{b\mu};
\ \ \ \widehat{C}^{\alpha}_{\mu c} := C^{\alpha}_{\mu c}, \ \ \ \widehat{C}^{a}_{bc}: =
\frac{1}{2}(C^{a}_{bc} + C^{a}_{cb}).\end{equation}
\end{description}

We shall denote the horizontal (vertical) covariant derivative of $\widetilde{D}$ and
$\widehat{D}$ by \linebreak \lq\lq\,\,$\widetilde{|}$\,\rq\rq
(\lq\lq\,\,$\widetilde{||}$\,\rq\rq) and \lq\lq\,\,$\widehat{|}$\,\rq\rq (\lq\lq\,\,$\widehat{||}$\,\rq\rq)
respectively.
\end{definition}

It follows immediately from the above definition that
 \begin{equation}\label{zero}{\lambda}^{\alpha}\!\,\,_{\widetilde{||}c} =
{\lambda}^{\alpha}\!\,\,_{\widehat{||}c} =
{\lambda}^{\alpha}\!\,_{||c} = 0, \ \ \ \ {\lambda}^{a}\!\,\,_{\widetilde{|}\mu} =
{\lambda}^{a}\!\,\,_{\widehat{|}\mu} =
{\lambda}^{a}\!\,_{|\mu} = 0;\end{equation}
\begin{equation}\label{half}\lambda^{\alpha}\!\,\,_{\widetilde{|}\mu} =
\lambda^{\beta}\Lambda^{\alpha}_{\mu\beta}, \ \ \ \ \lambda^{\alpha}\!\,\,_{\widehat{|}\mu} =
\frac{1}{2} \, \lambda^{\alpha}\!\,\,_{\widetilde{|}\mu};
\ \ \ \ \ \ \lambda^{a}\!\,\,_{\widetilde{||}c} = \lambda^{b}T^{a}_{cb},
\ \ \ \ \lambda^{a}\!\,\,_{\widehat{||}c} = \frac{1}{2} \, \lambda^{a}\!\,\,_{\widetilde{||}c}.\end{equation}

As easily checked, we also have

\begin{proposition} The covariant derivatives of the metric ${\bf G}$
with respect to the dual and symmetric $d$-connections $\widetilde{D}$ and $\widehat{D}$ are given respectively by:
\begin{equation}g_{\alpha\beta\widetilde{|}\mu} = \Lambda_{\alpha\beta\mu} + \Lambda_{\beta\alpha\mu}, \ \
g_{\alpha\beta\widetilde{||}c} = g_{ab\widetilde{|}\mu} = 0, \ \
g_{ab\widetilde{||}c} = T_{abc} + T_{bac};\end{equation}
\begin{equation}g_{\alpha\beta\widehat{|}\mu} = \frac{1}{2} \, g_{\alpha\beta\widetilde{|}\mu}
, \ \ \ \ \ \ \ g_{\alpha\beta\widehat{||}c} = g_{ab\widehat{|}\mu} = 0, \ \ \ \ \
g_{ab\widehat{||}c} = \frac{1}{2} \, g_{ab\widetilde{||}c}.\end{equation}
Consequently, $\widetilde{D}$ and $\widehat{D}$ are non-metric connections.
\end{proposition}

We end this section with the following tables.

\begin{center} {\bf Table 1: Fundamental connections of the EAP-space}\\[0.2cm]
\footnotesize{\begin{tabular}
{|c|c|c|c|c|c|c|c|c|c|c|c|}\hline
&&&&\\
{\bf Connection}&{\bf Coefficients}&{\bf Covariant}
 &{\bf Torsion}&{\bf Metricity}\\
 &&{\bf derivative}&&\\[0.2cm]\hline
&&&&\\
{\bf Natural}&\footnotesize${( \ \overcirc{\Gamma}^{\alpha}_{\mu\nu},
\ \overcirc{\Gamma}^{a}_{b\nu}, \ \overcirc{C}^{\alpha}_{\mu c}, \ \overcirc{C}^{a}_{bc})}
$&$\overcirc\atop{|}$ \ \ $\overcirc\atop||$&
\footnotesize${(0, R^{a}_{\mu\nu}, \ \overcirc{C}^{\alpha}_{\mu c},
\ \overcirc{P}^{a}_{\mu c}, 0)}$&metric\\[0.2 cm]\hline
&&&&\\
{\bf Canonical}&\footnotesize${(\Gamma^{\alpha}_{\mu\nu}, \Gamma^{a}_{b\nu},
C^{\alpha}_{\mu c}, C^{a}_{bc})}$&$|$ \ \ $||$
 &\footnotesize${(\Lambda^{\alpha}_{\mu\, \nu}, R^{a}_{\mu\nu},
C^{\alpha}_{\mu c}, P^{a}_{\mu c}, T^{a}_{bc})}$&metric\\[0.2cm]\hline
 &&&&\\
{\bf Dual}&\footnotesize${(\Gamma^{\alpha}_{\nu\mu}, \Gamma^{a}_{b\nu},
C^{\alpha}_{\mu c}, C^{a}_{cb})}$&$\begin{array}{cc}\widetilde {}\\[-0.3cm]|\end{array}$
$\begin{array}{cc}\widetilde {}\\[-0.3cm]||\end{array}$
&\footnotesize${(- \Lambda^{\alpha}_{\mu\nu},
R^{a}_{\mu\nu}, C^{\alpha}_{\mu c}, P^{a}_{\mu c}, - T^{a}_{bc})}$&non-metric\\[0.2cm]\hline
  &&&&\\
{\bf Symmetric}&\footnotesize${(\Gamma^{\alpha}_{(\mu\nu)}, \Gamma^{a}_{b\nu}, C^{\alpha}_{\mu c}, C^{a}_{(bc)})}$
&$\begin{array}{cc}
\widehat {}\\[-0.3cm]|\end{array}$$\begin{array}{cc}\widehat {}\\[-0.3cm]||\end{array}$
&\footnotesize${(0,
R^{a}_{\mu\nu}, C^{\alpha}_{\mu c}, P^{a}_{\mu c}, 0)}$&non-metric\\[0.2cm]\hline
\end{tabular}}
\end{center}
\par
\vspace{0.2 cm}

The next table gives a comparison between the classical AP-space and the EAP-space.
We shall refer to the Riemannian connection in the classical AP-space and the natural metric $d$-connection
in the EAP-space simply as the metric connection. Moreover, we consider only the metric and the canonical
connections in both spaces. We also set
$L^{a}_{b\nu}: = \frac{1}{2} \ g^{ac}(\delta_{\nu} g_{bc} - g_{dc} \ \dot{\partial_b}N^{d}_{\nu} -
g_{bd} \ \dot{\partial_c}N^{d}_{\nu}).$

\newpage

\begin{center} {\bf Table 2: Comparison between classical AP-geometry and EAP-geometry}\\[0.15 cm]
\footnotesize{\begin{tabular}
{|c|c|c|c|c|c|c|}\hline
&&\\
 \ \ &{\bf Classical AP-geometry}&{\bf EAP-geometry}
\\[0.3cm]\hline
&&\\
{\bf Underlying space}&$M$&$TM$
\\[0.2cm]\hline
&&\\
{\bf Building blocks}&${\lambda}^{\alpha}(x)$&$N^{a}_{\mu}(x, y), \ \ \, \lambda(x, y) =
(\lambda^{\alpha}(x, y), \ \lambda^{a}(x, y))$
\\[0.2cm]\hline
&&\\
{\bf Metric}&$g_{\mu\nu} = \ \undersym{\lambda}{i}_{\mu} \
\undersym{\lambda}{i}_{\nu}$ &${\bf G} = (g_{\mu\nu}, \ g_{ab});$
\\[0.01 cm]
&&\\
\ \ & \ \ & $g_{\mu\nu} = \
\undersym{\lambda}{i}_{\mu} \ \undersym{\lambda}{i}_{\nu}, \ \
g_{ab} = \ \undersym{\lambda}{i}_{a} \ \undersym{\lambda}{i}_{b}$
\\[0.2cm]\hline
&&\\
{\bf Metric}&\footnotesize${\overcirc{\Gamma}^{\alpha}_{\mu\nu} =
\frac{1}{2} \, g^{\alpha\epsilon}(\partial_{\mu}g_{\nu\epsilon} +
\partial_{\nu}g_{\mu\epsilon} - \partial_{\epsilon}g_{\mu\nu})}$&$\overcirc{D} =
( \ \overcirc{\Gamma}^{\alpha}_{\mu\nu}, \ \ \overcirc{\Gamma}^{a}_{b\nu}, \ \ \overcirc{C}^{\alpha}_{\mu c}, \
\ \overcirc{C}^{a}_{bc});$
\\[0.01 cm]
&&\\
{\bf connection}& \ \ & \footnotesize${\overcirc{\Gamma}^{\alpha}_{\mu\nu} =
\frac{1}{2} \, g^{\alpha\epsilon}(\delta_{\mu}g_{\nu\epsilon} +
\delta_{\nu}g_{\mu\epsilon} - \delta_{\epsilon}g_{\mu\nu})},$
\\[0.01 cm]
&&\\
\ \ & \ \ & \ $\overcirc{\Gamma}^{a}_{b\nu} = \dot{\partial_b}N^{a}_{\nu} + L^{a}_{b\nu}; \ \
\ \ \overcirc{C}^{\alpha}_{\mu c} = \frac{1}{2} \, g^{\alpha\epsilon}\dot{\partial_c}g_{\mu\epsilon},$

\\[0.01 cm]

&&\\
\ \ & \ \ & \footnotesize{$\overcirc{C}^{a}_{bc} =
\frac{1}{2} \, g^{ad}(\dot{\partial}_{b}g_{cd} +
\dot{\partial}_{c}g_{bd} -
\dot{\partial}_{d}g_{bc})$}
\\[0.2cm]\hline
&&\\
{\bf Canonical}&$\Gamma^{\alpha}_{\mu\nu} = \ \undersym{\lambda}{i}^{\alpha} (\partial_{\nu} \
\undersym{\lambda}{i}_{\mu})$& $D = (\Gamma^{\alpha}_{\mu\nu}, \ \Gamma^{a}_{b\nu}, \
C^{\alpha}_{\mu c}, \ C^{a}_{bc});$
\\[0.01 cm]
&&\\
{\bf connection}& \ \ &$\Gamma^{\alpha}_{\mu\nu} = \ \undersym{\lambda}{i}^{\alpha}(\delta_{\nu} \
\undersym{\lambda}{i}_{\mu}); \ \ \ \ \ \ \ \Gamma^{a}_{b\nu} = \ \undersym{\lambda}{i}^{a}(\delta_{\nu} \
\undersym{\lambda}{i}_{b}),$
\\[0.01 cm]
&&\\
\ \ & \ \ & $C^{\alpha}_{\mu c} = \
\undersym{\lambda}{i}^{\alpha}(\dot{\partial_{c}} \
\undersym{\lambda}{i}_{\mu}); \ \ \ \ \ \ \ C^{a}_{bc} = \
\undersym{\lambda}{i}^{a}(\dot{\partial_{c}} \
\undersym{\lambda}{i}_{b})$
\\[0.2 cm]\hline
&&\\
{\bf AP-condition}& ${\lambda}^{\alpha} \,\!_{|\mu} =
0$&${\lambda}^{\alpha}\, \!_{|\mu} = \lambda^{\alpha}\, \!_{||c} = 0,$
\\[0.01 cm]
&&\\
\ \ & \ \ & $\lambda^{a}\, \!_{|\mu} = \lambda^{a}\, \!_{||c} = 0$
\\[0.2 cm]\hline
&&\\
{\bf Torsion}& $\Lambda^{\alpha}_{\mu\nu} = \Gamma^{\alpha}_{\mu\nu} -
\Gamma^{\alpha}_{\nu\mu}$&${\bf T} = (\Lambda^{\alpha}_{\mu\nu}, \ R^{a}_{\mu\nu},
\ C^{\alpha}_{\mu c}, \ P^{a}_{\mu c}, \ T^{a}_{bc});$
\\[0.01 cm]
&&\\
\ \ & \ \ &$\Lambda^{\alpha}_{\mu\nu} = \Gamma^{\alpha}_{\mu\nu} - \Gamma^{\alpha}_{\nu\mu}; \
\ \, \, R^{a}_{\mu\nu} = \delta_{\nu}N^{a}_{\mu} - \delta_{\mu}N^{a}_{\nu},$
\\[0.01 cm]
&&\\
\ \ & \ \ &$P^{a}_{\mu c} = \dot{\partial_c}N^{a}_{\mu} - \Gamma^{a}_{c\mu}; \ \ \
T^{a}_{bc} = C^{a}_{bc} - C^{a}_{cb}$
\\[0.2 cm]\hline
&&\\
{\bf Contorsion}& $\gamma^{\alpha}_{\mu\nu} = \Gamma^{\alpha}_{\mu\nu} - \
\overcirc{\Gamma}^{\alpha}_{\mu\nu}$&${\bf C} = (\gamma^{\alpha}_{\mu\nu}, \ \gamma^{a}_{b\nu}, \
\gamma^{\alpha}_{\mu c}, \ \gamma^{a}_{bc});$
\\[0.01 cm]
&&\\
\ \ & \ \ &$\gamma^{\alpha}_{\mu\nu} =
\Gamma^{\alpha}_{\mu\nu} - \ \overcirc{\Gamma}^{\alpha}_{\nu\mu}; \ \ \ \ \gamma^{a}_{b\nu} =
\Gamma^{a}_{b\nu} - \ \overcirc{\Gamma}^{a}_{b\nu},$
\\[0.01 cm]
&&\\
\ \ & \ \ & $\gamma^{\alpha}_{\mu c} = C^{\alpha}_{\mu c} - \
\overcirc{C}^{\alpha}_{\mu c}; \ \ \ \ \ \gamma^{a}_{bc} = C^{a}_{bc} - \ \overcirc{C}^{a}_{bc}$
\\[0.2 cm]\hline
&&\\
{\bf Basic vector}& $C_{\mu} = \Lambda^{\alpha}_{\mu\alpha} =
\gamma^{\alpha}_{\mu\alpha}$& ${\bf B} = (C_{\mu}, \ C_{a});$
\\[0.01 cm]
&&\\
\ \ & \ \ & $C_{\mu} = \Lambda^{\alpha}_{\mu\alpha} = \gamma^{\alpha}_{\mu\alpha}; \ \ \ \
C_{a} = T^{d}_{ad} = \gamma^{d}_{ad}$
\\[0.2 cm]\hline
\end{tabular}}
\end{center}
\end{mysect}


\newpage

\begin{mysect}{Curvature tensors in the EAP-space}
Let $(TM, \ \lambda)$ be an EAP space. Let $D$ be the cannonical $d$-connection.
Let \ $\overcirc{D}$, $\widetilde{D}$ and
$\widehat{D}$ be the natural metric $d$-connection, the dual $d$-connection and the
symmetric $d$-connection respectively.
The curvature tensors of the four $d$-connections will be denoted respectively
by ${\bf R}$, \ $\overcirc{{\bf R}}$, ${\bf \widetilde{R}}$ and ${\bf \widehat{R}}$. In this section, we carry out the
task of calculating the curvature tensors together with their contractions.

\begin{lemma}\label{CFA} The following commutation formulae hold:
\begin{description}
\item [(a)] ${\lambda}^{\alpha}\!\, _{|\nu|\mu} - {\lambda}^{\alpha}\!\, _{|\mu|\nu} =
R^{\alpha}_{\beta\mu\nu}{\lambda}^{\beta} -
\Lambda^{\beta}_{\nu\mu}{\lambda}^{\alpha}\!\,_{|\beta} -
R^{d}_{\nu\mu}{\lambda}^{\alpha}\!\,_{||d}$
\item [(b)] ${\lambda}^{\alpha}\!\, _{|\nu||c} - {\lambda}^{\alpha}\!\, _{||c|\nu} =
P^{\alpha}_{\beta\nu c}{\lambda}^{\beta} -
C^{\beta}_{\nu c}{\lambda}^{\alpha}\!\, _{|\beta} - P^{d}_{\nu c}{\lambda}^{\alpha}\!\,_{||d}$
\item [(c)] ${\lambda}^{\alpha}\!\, _{||b||c} - {\lambda}^{\alpha}\!\, _{||c||b} =
S^{\alpha}_{\beta cb}{\lambda}^{\beta} -
T^{d}_{bc}{\lambda}^{\alpha}\!\, _{||d}$
\item [(d)] ${\lambda}^{a}\!\, _{|\nu|\mu} - {\lambda}^{a}\!\,_{|\mu|\nu} =
R^{a}_{b\mu\nu}{\lambda}^{b} -
\Lambda^{\beta}_{\nu\mu}{\lambda}^{a}\!\, _{|\beta} - R^{d}_{\nu\mu}{\lambda}^{a}\!\,_{||d}$
\item [(e)] \ ${\lambda}^{a}\!\,_{|\nu||c} - {\lambda}^{a}\!\,_{||c|\nu} =
P^{a}_{b\nu c}{\lambda}^{b} -
C^{\beta}_{\nu c}{\lambda}^{a}\!\,_{|\beta} - P^{d}_{\nu c}{\lambda}^{a}\!\, _{||d}$
\item [(f)] \ ${\lambda}^{a}\!\, _{||b||c} - {\lambda}^{a}\!\, _{||c||b} =
S^{a}_{dcb}{\lambda}^{d} -
T^{d}_{bc}{\lambda}^{a}\!\,_{||d},$
\end{description}
\end{lemma}
In view of (\ref{inverse}), Corollary \ref{vanish} and Theorem \ref{apc},
Lemma \ref{CFA} directly implies that

\begin{theorem}\label{v} The curvature tensor ${\bf {R}} =
({R}^{\alpha}_{\beta\mu\nu}, \!\, {R}^{a}_{b \mu\nu}, \!\, {P}^{\alpha}_{\beta\nu c}, \!\, {P}^{a}_{b \nu c},
\!\, {S}^{\alpha}_{\beta bc}, \!\, {S}^{a}_{bcd})$ of the \linebreak canonical $d$-connection vanishes identically.
\end{theorem}

\begin{corollary}\label{Bianchi} The following identities hold:

\begin{equation}\label{L}\Lambda^{\alpha}_{\beta\mu|\alpha} =
(C_{\beta|\mu} - C_{\mu|\beta}) + C_{\epsilon}\Lambda^{\epsilon}_{\beta\mu} + \
\mathfrak{S}_{\beta\mu\alpha} R^{a}_{\mu\beta} C^{\alpha}_{\alpha a}\end{equation}

\begin{equation}\label{T}T^{d}_{bc||d} = (C_{b||c} - C_{c||b}) + C_{d}T^{d}_{bc}.\end{equation}

\end{corollary}
\begin{proof} The Bianchi identities \cite{V} applied to the canonical $d$-connection $D$ give
\begin{equation}\label{FBI}\mathfrak{S}_{\beta, \mu, \nu} \ (\Lambda^{\alpha}_{\beta\mu|\nu} +
\Lambda^{\epsilon}_{\mu\nu}\Lambda^{\alpha}_{\beta\epsilon} + R^{a}_{\beta\mu}C^{\alpha}_{\nu a}) = 0\end{equation}
and
\begin{equation}\label{VFBI}\mathfrak{S}_{b, c, d} \ (T^{a}_{bc||d} +
T^{e}_{cd}T^{a}_{be}) = 0,\end{equation}
where the notation $\mathfrak{S}_{\beta, \mu, \nu}$ denotes a cyclic permutation
on the indices $\beta, \mu, \nu$ and summation. (\ref{L}) and (\ref{T}) are
obtained by setting $\alpha = \nu$ and $a = d$ in (\ref{FBI}) and
(\ref{VFBI}) respectively.
\end{proof}

By applying the comutation formula {\bf (c)} of Lemma \ref{CFA} with
respect to the dual and symmetric $d$-connections respectively, taking into account (\ref{inverse}) and (\ref{zero}), we obtain
$$\widetilde{S}^{\alpha}_{\beta bc} = \widehat{S}^{\alpha}_{\beta bc} = 0.$$

This could be also deduced from Theorem \ref{CUR} {\bf (e)} and Theorem \ref{v}, noting that $C^{\alpha}_{\mu c} =
\widetilde{C}^{\alpha}_{\mu c} = \widehat{C}^{\alpha}_{\mu c}$.

\newpage

In Theorems \ref{n}, \ref{d} and \ref{S} below, concerning the curvature tensors of the $d$-connections $\, \overcirc{D}$, $\widetilde{D}$ and $\widehat{D}$,
we will make use of Theorem \ref{v}, namely that the
curvature tensors of the canonical $d$-connection {\bf vanish} identically.

\begin{theorem} \label{n}The curvature tensors of the natural metric $d$-connection \ $\overcirc{D}$
can be expressed in the form:
\begin{description}
\item [(a)] $\overcirc{R}^{\alpha}_{\beta\mu\nu} = (\gamma^{\alpha}_{\beta\mu|\nu} -
\gamma^{\alpha}_{\beta\nu|\mu}) +
(\gamma^{\epsilon}_{\beta\nu}\gamma^{\alpha}_{\epsilon\mu} -
\gamma^{\epsilon}_{\beta\mu}\gamma^{\alpha}_{\epsilon\nu})
- \gamma^{\alpha}_{\beta\epsilon} \Lambda^{\epsilon}_{\nu\mu} -  \gamma^{\alpha}_{\beta d}R^{d}_{\nu\mu},$
\item [(b)] $\overcirc{R}^{a}_{b\mu\nu} = (\gamma^{a}_{b\mu|\nu} -
\gamma^{a}_{b\nu|\mu}) +
(\gamma^{d}_{b\nu}\gamma^{a}_{d\mu} - \gamma^{d}_{b\mu}\gamma^{a}_{d\nu}) -  \gamma^{a}_{b\epsilon}
\Lambda^{\epsilon}_{\nu\mu} - \gamma^{a}_{bd}R^{d}_{\nu\mu},$

\item [(c)] $\overcirc{P}^{\alpha}_{\beta\nu c} = (\gamma^{\alpha}_{\beta c|\nu} -
\gamma^{\alpha}_{\beta\nu||c}) +
(\gamma^{\epsilon}_{\beta\nu}\gamma^{\alpha}_{\epsilon c} -
\gamma^{\epsilon}_{\beta c}\gamma^{\alpha}_{\epsilon\nu}) - \gamma^{\alpha}_{\beta\epsilon}
C^{\epsilon}_{\nu c} - \gamma^{\alpha}_{\beta d}P^{d}_{\nu c},$
\item [(d)] $\overcirc{P}^{a}_{b\nu c} = (\gamma^{a}_{bc|\nu} -
\gamma^{a}_{b\nu||c}) +
(\gamma^{d}_{b\nu}\gamma^{a}_{dc} - \gamma^{d}_{bc}\gamma^{a}_{d\nu}) -
 \gamma^{a}_{b\epsilon} C^{\epsilon}_{\nu c} - \gamma^{a}_{bd}P^{d}_{\nu c},$

\item [(e)] $\overcirc{S}^{\alpha}_{\beta bc} = (\gamma^{\alpha}_{\beta b||c} -
\gamma^{\alpha}_{\beta c||b}) + (\gamma^{\epsilon}_{\beta c}\gamma^{\alpha}_{\epsilon b}
- \gamma^{\epsilon}_{\beta b}\gamma^{\alpha}_{\epsilon c}) - \gamma^{\alpha}_{\beta d}
T^{d}_{cb},$
\item [(f)] $\overcirc{S}^{a}_{bcd} = (\gamma^{a}_{bc||d} -
\gamma^{a}_{bd||c}) + (\gamma^{e}_{bd}\gamma^{a}_{ec} - \gamma^{e}_{bc}\gamma^{a}_{ed}) -  \gamma^{a}_{be}
T^{e}_{dc}.$
\end{description}

Consequently,
\begin{description}
\item [(g)] $\overcirc{R}_{\beta\mu}:= \ \overcirc{R}^{\alpha}_{\beta\mu\alpha} = (\gamma^{\alpha}_{\beta\mu|\alpha} -
C_{\beta|\mu}) - C_{\epsilon}\gamma^{\epsilon}_{\beta\mu} +
\gamma^{\alpha}_{\beta\epsilon}\gamma^{\epsilon}_{\mu\alpha} -
\gamma^{\alpha}_{\beta d}R^{d}_{\alpha\mu},$
\item [(h)] $\overcirc{\cal {R}}: = g^{\beta\mu} \ \overcirc{R}_{\beta\mu} = \frac{1}{2}(\Omega^{\alpha\mu}\!\,_{\mu|\alpha} -
C_{\alpha}\Omega^{\alpha\mu}\!\,_{\mu}) - C^{\mu}\!\,_{|\mu} + \gamma^{\alpha\mu}\!\,_{\epsilon}
\gamma^{\epsilon}_{\mu\alpha} - \gamma^{\alpha\mu}\!\,_{d}R^{d}_{\alpha\mu},$
\item [(i)] $\overcirc{P}_{\beta c}: = - \ \overcirc{P}^{\alpha}_{\beta\alpha c} = (C_{\beta||c} - \gamma^{\alpha}_{\beta c|\alpha})
 + C_{\epsilon}\gamma^{\epsilon}_{\beta c} +
\gamma^{\alpha}_{\beta\epsilon}(C^{\epsilon}_{\alpha c} - \gamma^{\epsilon}_{\alpha c})
+ \gamma^{\alpha}_{\beta d}P^{d}_{\alpha c},$
\item [(j)]  $\overcirc{P}_{b\nu}: = \ \overcirc{P}^{d}_{b\nu d} =
(C_{b|\nu} - \gamma^{d}_{b\nu||d}) + C_{d}\gamma^{d}_{b\nu} -
\gamma^{d}_{be}\gamma^{e}_{d\nu} - \gamma^{d}_{b\epsilon}C^{\epsilon}_{\nu d}
- \gamma^{e}_{bd}P^{d}_{\nu e},$
\item [(k)] $\overcirc{S}_{bc}: = \ \overcirc{S}^{d}_{bcd} = (\gamma^{d}_{bc||d} - C_{b||c}) - C_{d}\gamma^{d}_{bc} +
\gamma^{d}_{be}\gamma^{e}_{cd},$
\item [(l)] $\overcirc{\cal{S}} := g^{bc} \ \overcirc{S}_{bc} = \frac{1}{2}(\Omega^{ad}\!\,_{d||a} -
C_{a}\Omega^{ad}\!\,_{d}) - C^{d}\!\,_{||d} + \gamma^{ad}\!\,_{c}
\gamma^{c}_{da},$
\end{description}
where \ $\Omega^{\alpha}_{\beta\mu}: = \gamma^{\alpha}_{\beta\mu} + \gamma^{\alpha}_{\mu\beta}, \ \
\Omega^{a}_{bc}: = \gamma^{a}_{bc} + \gamma^{a}_{cb}.$

\end{theorem}

\begin{proof} We prove {\bf (a)} and {\bf (c)} only. The other formulae of the first part
are proved in a similar manner. The second part is obtained directly by applying the suitable contractions.
\begin{description}
\item [(a)] We have
\begin{eqnarray*} \overcirc{R}^{\alpha}_{\beta\mu\nu} &=&
\delta_{\mu} \ \overcirc{\Gamma}^{\alpha}_{\beta\nu} -
\delta_{\nu} \ \overcirc{\Gamma}^{\alpha}_
{\beta\mu} + \ \overcirc{\Gamma}^{\epsilon}_{\beta\nu}
\ \overcirc{\Gamma}^{\alpha}_{\epsilon\mu} - \ \overcirc{\Gamma}^{\epsilon}_{\beta\mu}
\ \overcirc{\Gamma}^{\alpha}_{\epsilon\nu} + \ \overcirc{C}^{\alpha}_{\beta d}R^{d}_{\nu\mu}\\
&=& \delta_{\mu}({\Gamma}^{\alpha}_{\beta\nu} - \gamma^{\alpha}_{\beta\nu}) -
\delta_{\nu} ({\Gamma}^{\alpha}_{\beta\mu} - \gamma^{\alpha}_{\beta\mu})
+ ({\Gamma}^{\epsilon}_{\beta\nu} - \gamma^{\epsilon}_{\beta\nu})
({\Gamma}^{\alpha}_{\epsilon\mu} - \gamma^{\alpha}_{\epsilon\mu}) -
({\Gamma}^{\epsilon}_{\beta\mu} - \gamma^{\epsilon}_{\beta\mu})\\
&& \ (\Gamma^{\alpha}_{\epsilon\nu} - \gamma^{\alpha}_{\epsilon\nu}) +
({C}^{\alpha}_{\beta d} - \gamma^{\alpha}_{\beta d})R^{d}_{\nu\mu}\\
&=& R^{\alpha}_{\beta\mu\nu} + (\delta_{\nu}\gamma^{\alpha}_{\beta\mu} +
\gamma^{\epsilon}_{\beta\mu}\Gamma^{\alpha}_{\epsilon\nu} -
\gamma^{\alpha}_{\epsilon\mu}\Gamma^{\epsilon}_{\beta\nu}) -
(\delta_{\mu}\gamma^{\alpha}_{\beta\nu} +
\gamma^{\epsilon}_{\beta\nu}\Gamma^{\alpha}_{\epsilon\mu} -
\gamma^{\alpha}_{\epsilon\nu}\Gamma^{\epsilon}_{\beta\mu})\\
&& - \ \gamma^{\alpha}_{\beta d} R^{d}_{\nu\mu} + (\gamma^{\epsilon}_{\beta\nu}
\gamma^{\alpha}_{\epsilon\mu} - \gamma^{\epsilon}_{\beta\mu}
\gamma^{\alpha}_{\epsilon\nu})\\
&=& (\gamma^{\alpha}_{\beta\mu|\nu} - \gamma^{\alpha}_{\beta\nu|\mu}) +
(\gamma^{\epsilon}_{\beta\nu}\gamma^{\alpha}_{\epsilon\mu} - \gamma^{\epsilon}_{\beta\mu}
\gamma^{\alpha}_{\epsilon\nu}) - \gamma^{\alpha}_{\beta\epsilon}\Lambda^{\epsilon}_{\nu\mu}
- \gamma^{\alpha}_{\beta d} R^{d}_{\nu\mu}\\
\end{eqnarray*}
\item [(c)] We have
\begin{eqnarray*} \overcirc{P}^{\alpha}_{\beta\nu c} &=& \dot{\partial_c}
 \ \overcirc{\Gamma}^{\alpha}_{\beta\nu} - \ \overcirc{C}^{\alpha}_{\beta c{o\atop{|}}\nu} +
\ \overcirc{C}^{\alpha}_{\beta d} \ \overcirc{P}^{d}_{\nu c}\\
&=& \dot{\partial_c}\Gamma^{\alpha}_{\beta\nu} - \dot{\partial_c}
\gamma^{\alpha}_{\beta\nu} - C^{\alpha}_{\beta c|\nu} + (C^{\alpha}_{\beta c|\nu} -
\ \overcirc{C}^{\alpha}_{\beta c{o\atop{|}}\nu}) + ({C}^{\alpha}_{\beta d} - \gamma^{\alpha}_{\beta d})(P^{d}_{\nu c} + \gamma^{d}_{c\nu})\\
&=& P^{\alpha}_{\beta\nu c} - \dot{\partial_c}\gamma^{\alpha}_{\beta\nu}
+ \{C^{\alpha}_{\beta c|\nu}- ({C}^{\alpha}_{\beta c} -
\gamma^{\alpha}_{\beta c})_{{o\atop{|}}\nu}\}
+ C^{\alpha}_{\beta d}\gamma^{d}_{c\nu} -
\gamma^{\alpha}_{\beta d}\gamma^{d}_{c \nu} - \gamma^{\alpha}_{\beta d}{P}^{d}_{\nu c}\\
&=& - \ \dot{\partial_c}\gamma^{\alpha}_{\beta\nu}
+ (C^{\epsilon}_{\beta c}\gamma^{\alpha}_{\epsilon \nu} - {C}^{\alpha}_{\epsilon c}
\gamma^{\epsilon}_{\beta \nu} - C^{\alpha}_{\beta d}\gamma^{d}_{c\nu}) +
(\gamma^{\alpha}_{\beta c|\nu} - \gamma^{\epsilon}_{\beta c}
\gamma^{\alpha}_{\epsilon\nu} +  \gamma^{\alpha}_{\epsilon c}\gamma^{\epsilon}_{\beta\nu}\\
&& + \ \gamma^{\alpha}_{\beta d}\gamma^{d}_{c \nu}) + C^{\alpha}_{\beta d} \gamma^{d}_{c\nu} - \gamma^{\alpha}_{\beta d}\gamma^{d}_{c \nu}
- \gamma^{\alpha}_{\beta d}{P}^{d}_{\nu c}\\ &=& \gamma^{\alpha}_{\beta c|\nu} - ( \dot{\partial_c}\gamma^{\alpha}_{\beta\nu} +
\gamma^{\epsilon}_{\beta\nu}C^{\alpha}_{\epsilon c} -
\gamma^{\alpha}_{\epsilon \nu}{C}^{\epsilon}_{\beta c}) +  (\gamma^{\epsilon}_{\beta\nu}
\gamma^{\alpha}_{\epsilon c} - \gamma^{\epsilon}_{\beta c}\gamma^{\alpha}_{\epsilon\nu}) - \gamma^{\alpha}_{\beta d} P^{d}_{\nu c}\\
&=& (\gamma^{\alpha}_{\beta c|\nu} - \gamma^{\alpha}_{\beta\nu||c}) + (\gamma^{\epsilon}_{\beta\nu}
\gamma^{\alpha}_{\epsilon c} - \gamma^{\epsilon}_{\beta c}\gamma^{\alpha}_{\epsilon\nu}) -
\gamma^{\alpha}_{\beta\epsilon}C^{\epsilon}_{\nu c} - \gamma^{\alpha}_{\beta d} P^{d}_{\nu c}\\
\end{eqnarray*}
\\[- 2.5 cm]\end{description}
\end{proof}

\begin{theorem}\label{d} The non-vanishing curvature tensors of the dual $d$-connection
$\widetilde{D}$ can be expressed  in the form:
\begin{description}
\item [(a)] $\widetilde{R}^{\alpha}_{\beta\nu\mu} = \Lambda^{\alpha}_{\mu\nu|\beta} +
\mathfrak{S}_{\beta, \mu, \nu} C^{\alpha}_{\beta a}R^{a}_{\mu\nu},$
\item [(b)]  $\widetilde{R}^{a}_{b\nu\mu} = R^{d}_{\mu\nu}T^{a}_{db},$
\item [(c)] $\widetilde{P}^{\alpha}_{\beta\mu c} = \Lambda^{\alpha}_{\mu\beta||c} +
\Lambda^{\alpha}_{\epsilon\beta} C^{\epsilon}_{\mu c},$
\item [(d)] $\widetilde{P}^{a}_{b\mu c} = T^{a}_{bc|\mu} + T^{a}_{db}P^{d}_{\mu c},$
\item [(e)] $\widetilde{S}^{a}_{bcd} = T^{a}_{dc||b}.$
\end{description}
Consequently,
\begin{description}
\item [(f)] $\widetilde{R}_{\beta\nu} := \widetilde{R}^{\alpha}_{\beta\nu\alpha} = - \ C_{\nu|\beta} +
\mathfrak{S}_{\beta, \nu, \alpha} C^{\alpha}_{\beta a}R^{a}_{\alpha\nu},$
\item [(g)] $\widetilde{\cal{R}} := g^{\beta\mu}\widetilde{R}_{\beta\mu} =  - \ C^{\mu}\!\,_{|\mu},$
\item [(h)] $\widetilde{P}_{\beta c} := - \widetilde{P}^{\alpha}_{\beta\alpha c} =
C_{\beta||c} + \Lambda^{\alpha}_{\beta\epsilon}C^{\epsilon}_{\alpha c},$
\item [(i)] $\widetilde{P}_{b\mu} :=  \widetilde{P}^{a}_{\beta\mu a} =
C_{b|\mu} + T^{a}_{db}P^{d}_{\mu a},$
\item [(j)]$\widetilde{S}_{bd} := \widetilde{S}^{a}_{bda} = - \ C_{d||b},$
\item [(k)] $\widetilde{\cal{S}} := g^{bd}\widetilde{S}_{bd} = - \ C^{d}\!\,_{||d}.$
\end{description}
\end{theorem}

\begin{proof} {\bf (b)} is a consequence of the commutation formula {\bf (d)} of Lemma \ref{CFA} applied to the dual
$d$-connection, taking into account (\ref{inverse}), (\ref{zero}) and (\ref{half}). {\bf (b)} could be also obtained from
Theorem \ref{CUR} {\bf (b)} and Theorem \ref{v}, noting that
$\Gamma^{a}_{b\mu} = \widetilde{\Gamma}^{a}_{b\mu}$ and
$C^{a}_{bd} = T^{a}_{bd}\, + \, \widetilde{C}^{a}_{bd}.$
We next prove {\bf (a)} and {\bf (c)} of the first part.
The second part follows immediately by applying the suitable contractions.
\begin{description}
\item[(a)] We have
\begin{eqnarray*} \tilde{R}^{\alpha}_{\beta\nu\mu} &=&
\delta_{\nu}\widetilde{\Gamma}^{\alpha}_{\beta\mu} -
\delta_{\mu}\widetilde{\Gamma}^{\alpha}_{\beta\nu} +
\widetilde{\Gamma}^{\epsilon}_{\beta\mu}\widetilde{\Gamma}^{\alpha}_{\epsilon\nu} -
\widetilde{\Gamma}^{\epsilon}_{\beta\nu}\widetilde{\Gamma}^{\alpha}_{\epsilon\mu} +
\widetilde{C}^{\alpha}_{\beta a}R^{a}_{\mu\nu}\\
&=& \delta_{\nu}{\Gamma}^{\alpha}_{\mu\beta} -
\delta_{\mu}{\Gamma}^{\alpha}_{\nu\beta} +
{\Gamma}^{\epsilon}_{\mu\beta}{\Gamma}^{\alpha}_{\nu\epsilon} -
{\Gamma}^{\epsilon}_{\nu\beta}{\Gamma}^{\alpha}_{\mu\epsilon} +
{C}^{\alpha}_{\beta a}R^{a}_{\mu\nu}\\
&=& \{\delta_{\nu}\Gamma^{\alpha}_{\mu\beta} +
\Gamma^{\epsilon}_{\mu\beta}(\Lambda^{\alpha}_{\nu\epsilon} +
\Gamma^{\alpha}_{\epsilon\nu})\} -  \{\delta_{\mu}\Gamma^{\alpha}_{\nu\beta} +
\Gamma^{\epsilon}_{\nu\beta}(\Lambda^{\alpha}_{\mu\epsilon} +
\ \Gamma^{\alpha}_{\epsilon\mu})\} + C^{\alpha}_{\beta a}R^{a}_{\mu\nu}\\
&=& (R^{\alpha}_{\mu\nu\beta} - C^{\alpha}_{\mu a} R^{a}_{\beta\nu} +
\delta_{\beta}\Gamma^{\alpha}_{\mu\nu} + \Gamma^{\epsilon}_{\mu\nu}
\Gamma^{\alpha}_{\epsilon\beta}) - (R^{\alpha}_{\nu\mu\beta} - C^{\alpha}_{\nu a}
R^{a}_{\beta\mu} + \delta_{\beta}\Gamma^{\alpha}_{\nu\mu}\\
&& + \ \Gamma^{\epsilon}_{\nu\mu}\Gamma^{\alpha}_{\epsilon\beta})
- \ (\Gamma^{\epsilon}_{\mu\beta}\Lambda^{\alpha}_{\epsilon\nu} +
\Gamma^{\epsilon}_{\nu\beta}\Lambda^{\alpha}_{\mu\epsilon}) +
C^{\alpha}_{\beta a}R^{a}_{\mu\nu}\\
&=& (\delta_{\beta}\Lambda^{\alpha}_{\mu\nu} + \Gamma^{\alpha}_{\epsilon\beta}
\Lambda^{\epsilon}_{\mu\nu} - \Gamma^{\epsilon}_{\mu\beta}
\Lambda^{\alpha}_{\epsilon\nu} - \Gamma^{\epsilon}_{\nu\beta}\Lambda^{\alpha}_
{\mu\epsilon}) + \mathfrak{S}_{\beta, \mu, \nu}C^{\alpha}_{\beta a}R^{a}_{\mu\nu}\\
&=& \Lambda^{\alpha}_{\mu\nu|\beta} +  \mathfrak{S}_{\beta, \mu, \nu}
C^{\alpha}_{\beta a}R^{a}_{\mu\nu}\\
\end{eqnarray*}
\\[- 2 cm] \item [(c)] We have
\begin{eqnarray*} \widetilde{P}^{\alpha}_{\beta\mu c} &=&
\dot{\partial_c}\widetilde{\Gamma}^{\alpha}_{\beta\mu} -
\widetilde{C}^{\alpha}_{\beta c\widetilde{|}\mu} +
\widetilde{C}^{\alpha}_{\beta d} \widetilde{P}^{d}_{\mu c} =
\dot{\partial_c}\Gamma^{\alpha}_{\mu\beta} - C^{\alpha}_{\beta c\widetilde{|}\mu} +
C^{\alpha}_{\beta d}P^{d}_{\mu c}\\
&=& (\dot{\partial_c}\Gamma^{\alpha}_{\beta\mu} - C^{\alpha}_{\beta c|\mu} +
C^{\alpha}_{\beta d}P^{d}_{\mu c}) +
\dot{\partial_c}\Lambda^{\alpha}_{\mu\beta} + (C^{\alpha}_{\beta c|\mu} -
{C}^{\alpha}_{\beta c\widetilde{|}\mu})\\
&=& P^{\alpha}_{\beta\mu c} + (\dot{\partial_c}\Lambda^{\alpha}_{\mu\beta} +
\Lambda^{\epsilon}_{\mu\beta}C^{\alpha}_{\epsilon c} - \Lambda^{\alpha}_{\mu\epsilon}
C^{\epsilon}_{\beta c})\\
&=& \Lambda^{\alpha}_{\mu\beta||c} + \Lambda^{\alpha}_{\epsilon\beta} C^{\epsilon}_{\mu c}\\
\end{eqnarray*}
\\[- 2.5 cm]\end{description}
\vspace{- 0.4cm}\end{proof}

\begin{theorem} \label{S}The non-vanishing curvature tensors of the symmetic $d$-connection
$\widehat{D}$ can be expressed
in the form
\begin{description}
\item [(a)] $\widehat{R}^{\alpha}_{\beta\nu\mu} =
\frac{1}{2}\, (\Lambda^{\alpha}_{\beta\nu|\mu} - \Lambda^{\alpha}_{\beta\mu|\nu}) +
\frac{1}{4}\, (\Lambda^{\epsilon}_{\beta\nu} \Lambda^{\alpha}_{\mu\epsilon} -
\Lambda^{\epsilon}_{\beta\mu} \Lambda^{\alpha}_{\nu\epsilon}) +
\frac{1}{2}\, (\Lambda^{\epsilon}_{\nu\mu} \Lambda^{\alpha}_{\beta\epsilon}),$
\item [(b)] $\widehat{R}^{a}_{b\nu\mu} = \frac{1}{2} \, \widetilde{R}^{a}_{b\nu\mu},$

\item [(c)] $\widehat{P}^{\alpha}_{\beta\mu c} = \frac{1}{2}\, \widetilde{P}^{\alpha}_{\beta\mu c},$

\item [(d)] $\widehat{P}^{a}_{b\mu c} = \frac{1}{2}\, \widetilde{P}^{a}_{b\mu c},$

\item [(e)] $\widehat{S}^{a}_{bcd} = \frac{1}{2}\, (T^{a}_{bc||d} - T^{a}_{bd||c}) +
\frac{1}{4}\, (T^{e}_{bc}T^{a}_{de} - T^{e}_{bd}T^{a}_{ce}) +
\frac{1}{2}\, (T^{e}_{dc}T^{a}_{eb}).$
\end{description}

Consequently,
\begin{description}
\item [(f)] $\widehat{R}_{\beta\nu} := \widehat{R}^{\alpha}_{\beta\nu\alpha} =
\frac{1}{2}\, \widetilde{R}_{\beta\nu} - \frac{1}{4}\, (C_{\alpha}
\Lambda^{\alpha}_{\nu\beta} + \Lambda^{\alpha}_{\nu\epsilon}\Lambda^{\epsilon}_{\alpha\beta}),$
\item [(g)] $\widehat{\cal{R}} := g^{\beta\nu}\widetilde{R}_{\beta\nu} = \frac{1}{2}\, \widetilde{\cal{R}} -
\frac{1}{4}\, \Lambda^{\alpha\beta}\!\,_{\epsilon}\Lambda^{\epsilon}_{\alpha\beta},$
\item [(h)] $\widehat{P}_{\beta c} := - \widehat{P}^{\alpha}_{\beta\alpha c} =
\frac{1}{2}\, \widetilde{P}_{\beta c},$
\item [(i)] $\widehat{P}_{b\mu} :=  \widehat{P}^{a}_{\beta\mu a} =
\frac{1}{2}\, \widetilde{P}_{b \mu},$
\item [(j)]$\widehat{S}_{bd} := \widehat{S}^{a}_{bda} = \frac{1}{2}\, \widetilde{S}_{bd} - \frac{1}{4}\,
(C_{a}T^{a}_{db} + T^{a}_{de}T^{e}_{ab}),$
\item [(k)] $\widehat{\cal{S}} := g^{bd}\widehat{S}_{bd} = \frac{1}{2}\, \widetilde{\cal{S}} - \frac{1}{4}\,
T^{ab}\!\,_{e}T^{e}_{ab}.$
\end{description}
\end{theorem}
\begin{proof} Similar to the proof of Theorems \ref{n} and \ref {d}.
\vspace{- 0.4 cm}\end{proof}
\end{mysect}


\begin{mysect}{Wanas tensors (W-tensors)}

The Wanas tensor, or simply the $W$-tensor, in the classical
AP-geometry, is a tensor which measures the non-commutativity of
covariant differentiations of the parallelization vector fields
$\undersym{\lambda}{i}$ with respect to the dual connection:
\begin{equation}\label{WT}W^{\alpha}_{\beta\nu\mu} : = \ \undersym{\lambda}{i}_{\beta}
( \ \undersym{\lambda}{i}^{\alpha}\!\,\,_{\widetilde{|}\nu\widetilde{|}\mu} -  \
\undersym{\lambda}{i}^{\alpha}\!\,\,_{\widetilde{|}\mu\widetilde{|}\nu})\end{equation}
This tensor explicitely contains the curvature and torsion tensors.
The $W$-tensor was first defined by M. Wanas \cite{aa}
and
has been used by F. Mikhail and M. Wanas \cite{aaa}
to construct a geometric theory unifying gravity and electromagnetism (GFT: generalized field theory).
The scalar Lagrangian function of the GFT is obtained by double contractions of the
tensor $W^{\alpha\beta}\!\,_{\nu\mu}$. The symmetric part of the field equation obtained
contains a second-order tensor representing the material distribution. This tensor is a pure geometric,
not a phenomological, object. The skew part of the field equation gives rise to
Maxwell-like equations. The use of the $W$-tensor has thus aided to construct
a geometric theory via {\it one single geometric entity}, which Einstein was seeking for \cite{E}.
Various significant applications (e.g \cite{aa}, \cite{MIW}, \cite{W}) have supported such a theory.
Recently, the authors of this paper investigated the most important properties of
this tensor in the context of classical AP-geometry \cite{AMR}. The $W$-tensor
was also studied by the present authors in the context of generalized Lagrange spaces \cite{WNA}.
\par
It should be noted that the $W$-tensor can be defined only in the context of AP-geometry and its generalized versions (cf. e.g \cite{WNA}, \cite{AMR}),
since it is defined
only in terms of the vector fields ${\lambda}$'s.
\par
Due to its importance in physical applications, we are going to investigate the properties of the $W$-tensor
in the present section. The $W$-tensor (\ref{WT}) can be generalized in the context of the EAP-geometry as follows.

\begin{definition} Let $(TM, \ \lambda)$ be an EAP-space.
For a given $d$-connection
$D = (\Gamma^{\alpha}_{\mu\nu}, \ \Gamma^{a}_{b\mu}, \ \linebreak C^{\alpha}_{\mu c}, \ C^{a}_{bc})$,
the $W$-tensor is given by
$$W = (W^{\alpha}_{\beta\nu\mu}, \ W^{a}_{b\nu\mu}, \ W^{\alpha}_{\beta\nu c}, \
W^{a}_{b\nu c}, \ W^{\alpha}_{\beta bc}, \
W^{a}_{bcd}),$$

\begin{description}
\item [(a)] the $hhh$-tensor $W^{\alpha}_{\beta\nu\mu}$ is defined by
the formula
$${\lambda}^{\alpha}\!\,_{|\nu|\mu} - {\lambda}^{\alpha}\!\,_{|\mu|\nu} =
{\lambda}^{\epsilon}W^{\alpha}_{\epsilon\nu\mu},$$
\item [(b)] the $hhv$-tensor $W^{a}_{b\nu\mu}$ is
defined by the formula
$${\lambda}^{a}\!\,_{|\nu|\mu} - {\lambda}^{a}\!\,_{|\mu|\nu} =
{\lambda}^{d}W^{a}_{d\nu\mu},$$
\item [(c)] the $vhh$-tensor $W^{\alpha}_{\beta\nu c}$ is defined by
the formula
$${\lambda}^{\alpha}\!\,_{|\nu||c} - {\lambda}^{\alpha}\!\,_{||c|\nu} =
{\lambda}^{\epsilon}W^{\alpha}_{\epsilon\nu c},$$
\item [(d)] the $vhv$-tensor $W^{a}_{b\nu c}$ is
defined by the formula
$${\lambda}^{a}\!\,_{|\nu||c} - {\lambda}^{a}\!\,_{||c|\nu} =
{\lambda}^{d}W^{a}_{d\nu c},$$

\item [(e)] the $vvh$-tensor $W^{\alpha}_{\beta bc}$ is defined by
the formula
$${\lambda}^{\alpha}\!\,_{||b||c} - {\lambda}^{\alpha}\!\,_{||c||b} =
{\lambda}^{\epsilon}W^{\alpha}_{\epsilon bc},$$

\item [(f)] the $vvv$-tensor $W^{a}_{bcd}$ is
defined by the formula
$${\lambda}^{a}\!\,_{||c||d} - {\lambda}^{a}\!\,_{||d||c} =
{\lambda}^{e}W^{a}_{ecd},$$
\end{description}
where \lq\lq\ $|$\rq\rq\, and \lq\lq \,$||$\rq\rq\, are the h- and v-covariant
derivatives with respect to the given $d$-connection $D$.
\end{definition}

Theorem \ref{apc}, together with (\ref{inverse}), directly implies that the
$W$-tensors of the canonical $d$-connection vanish identically.

\bigskip

In view of Theorem \ref{n}, we obtain

\begin{theorem} \label{nw}The $W$-tensors corresponding to the natural metric $d$-connection \ $\overcirc{D}$ are given by:
\begin{description}
\item [(a)] $\overcirc{W}^{\alpha}_{\beta\nu\mu} = (\gamma^{\alpha}_{\beta\mu|\nu} -
\gamma^{\alpha}_{\beta\nu|\mu}) +
(\gamma^{\epsilon}_{\beta\nu}\gamma^{\alpha}_{\epsilon\mu} -
\gamma^{\epsilon}_{\beta\mu}\gamma^{\alpha}_{\epsilon\nu}) - \gamma^{\alpha}_{\beta\epsilon}
\Lambda^{\epsilon}_{\nu\mu},$
\item [(b)] $\overcirc{W}^{a}_{b\nu\mu} = (\gamma^{a}_{b\mu|\nu} -
\gamma^{a}_{b\nu|\mu}) +
(\gamma^{d}_{b\nu}\gamma^{a}_{d\mu} - \gamma^{d}_{b\mu}\gamma^{a}_{d\nu}) - \gamma^{a}_{b\epsilon}
\Lambda^{\epsilon}_{\nu\mu},$

\item [(c)] $\overcirc{W}^{\alpha}_{\beta\nu c} = (\gamma^{\alpha}_{\beta c|\nu} -
\gamma^{\alpha}_{\beta\nu||c})  + (\gamma^{\epsilon}_{\beta \nu}\gamma^{\alpha}_{\epsilon c} -
\gamma^{\epsilon}_{\beta c}\gamma^{\alpha}_{\epsilon \nu}) +
(\gamma^{d}_{c\nu}\gamma^{\alpha}_{\beta d} - \gamma^{\epsilon}_{\nu c}\gamma^{\alpha}_{\beta\epsilon}),$

\item [(d)] $\overcirc{W}^{a}_{b\nu c} = (\gamma^{a}_{bc|\nu} -
\gamma^{a}_{b\nu||c}) +
(\gamma^{d}_{b\nu}\gamma^{a}_{dc} - \gamma^{d}_{bc}\gamma^{a}_{d\nu}) +
(\gamma^{d}_{c\nu}\gamma^{a}_{bd} - \gamma^{\epsilon}_{\nu c}\gamma^{a}_{b\epsilon}),$

\item [(e)] $\overcirc{W}^{\alpha}_{\beta bc} = \ \overcirc{S}^{\alpha}_{\beta cb},$
\item [(f)] $\overcirc{W}^{a}_{bcd} = \ \overcirc{S}^{a}_{bdc}.$
\end{description}
\end{theorem}

\begin{proof} We prove {\bf (c)} only. The other formulae are derived in a similar manner.
By definition, we have
$$\lambda^{\epsilon} \ \overcirc{W}^{\alpha}_{\epsilon\nu c} =
\lambda^{\epsilon} \ \overcirc{P}^{\alpha}_{\epsilon\nu c} -  \
\overcirc{C}^{\epsilon}_{\nu c}\lambda^{\alpha}\!\,_{{o\atop{|}}\epsilon} - \ \overcirc{P}^{d}_{\nu c}
\lambda^{\alpha}\!\,_{{o\atop{||}}d}.$$
\\[- 0.4 cm]Consequently, by (\ref{inverse}), (\ref{contortion}) and (\ref{torsion}), we obtain
\begin{eqnarray*} \ \overcirc{W}^{\alpha}_{\beta\nu c} &=&  \ \overcirc{P}^{\alpha}_{\beta\nu c} -
( \ \undersym{\lambda}{i}_{\beta} \ \undersym{\lambda}{i}^{\alpha}\!\,_{{o\atop{|}}\epsilon})
(C^{\epsilon}_{\nu c} - \gamma^{\epsilon}_{\nu c})
- ( \ \undersym{\lambda}{i}_{\beta} \ \undersym{\lambda}{i}^{\alpha}\!\,_{{o\atop{||}}d})(P^{d}_{\nu c} +
\gamma^{d}_{c\nu})\\
&=& \ \overcirc{P}^{\alpha}_{\beta\nu c} + C^{\epsilon}_{\nu c}\gamma^{\alpha}_{\beta\epsilon} +
P^{d}_{\nu c}\gamma^{\alpha}_{\beta d} -
(\gamma^{\epsilon}_{\nu c}\gamma^{\alpha}_{\beta\epsilon} - \gamma^{d}_{c\nu}\gamma^{\alpha}_{\beta d})\\
&=&   (\gamma^{\alpha}_{\beta c|\nu} -
\gamma^{\alpha}_{\beta\nu||c})  + (\gamma^{\epsilon}_{\beta\nu}\gamma^{\alpha}_{\epsilon c} -
\gamma^{\epsilon}_{\beta c}\gamma^{\alpha}_{\epsilon\nu}) +
(\gamma^{d}_{c\nu}\gamma^{\alpha}_{\beta d} - \gamma^{\epsilon}_{\nu c}\gamma^{\alpha}_{\beta\epsilon})\\
\\[- 2 cm]\end{eqnarray*}
\end{proof}

\bigskip

Since $\lambda^{a}\!\,\,_{\widetilde{|}\mu} = \lambda^{\alpha}\!\,\,_{\widetilde{||}c} = \lambda^{a}\!\,\,_{\widehat{|}\mu} =
\lambda^{\alpha}\!\,\,_{\widehat{||}c} = 0$,
it follows, by definition, that
$$\widetilde{W}^{a}_{b\nu\mu} = \widetilde{W}^{\alpha}_{\beta bc} = \widehat{W}^{a}_{b\nu\mu} =
\widehat{W}^{\alpha}_{\beta bc} = 0.$$
Proceeding as in Theorem \ref{nw}, taking into account Theorem \ref{d} and Theorem \ref{S}, we have the following

\begin{theorem}\label{dw}The non-vanishing $W$-tensors corresponding to the
dual $d$-connection $\widetilde{D}$ are given by:
\begin{description}
\item [(a)] $\widetilde{W}^{\alpha}_{\beta\nu\mu} = \Lambda^{\alpha}_{\nu\mu|\beta} +
\Lambda^{\epsilon}_{\nu\mu}\Lambda^{\alpha}_{\beta\epsilon} + \mathfrak{S}_{\nu, \mu, \beta}
C^{\alpha}_{\beta a}R^{a}_{\nu\mu},$
\item [(b)] $\widetilde{W}^{\alpha}_{\beta\nu c} = \Lambda^{\alpha}_{\nu\beta||c},$
\item [(c)] $\widetilde{W}^{a}_{b\nu c} = T^{a}_{bc|\nu},$
\item [(d)] $\widetilde{W}^{a}_{bdc} = T^{a}_{dc||b} + T^{e}_{dc}T^{a}_{be}.$
\end{description}
\end{theorem}

\begin{theorem}\label{Sw}The non-vanishing $W$-tensors corresponding to the
symmetric $d$-\linebreak connection $\widehat{D}$ are given by:
\begin{description}
\item [(a)] $\widehat{W}^{\alpha}_{\beta\nu\mu} = \widehat{R}^{\alpha}_{\beta\mu\nu},$
\item [(b)] $\widehat{W}^{\alpha}_{\beta\nu c} = \frac{1}{2} \ \widetilde
{W}^{\alpha}_{\beta\nu c},$
\item [(c)] $\widehat{W}^{a}_{b\nu c} = \frac{1}{2} \ \widetilde{W}^{a}_{b\nu c},$
\item [(d)] $\widehat{W}^{a}_{bcd} = \widehat{S}^{a}_{bdc}.$
\end{description}
\end{theorem}

It is clear by the above theorem that the $W$-tensors corresponding to the
symmetric $d$-connection give no new $d$-tensor fields.

\begin{remark} \em{The $W$-tensors corresponding to a given $d$-connection can be also
defined covariantly in the form
$$\lambda_{\beta|\mu|\nu} - \lambda_{\beta|\nu|\mu} = \lambda_{\epsilon} W^{\epsilon}_{\beta\mu\nu},$$
with similar expressions for the other counterparts. These expressions give the same formulae (up to a sign)
for the $W$-tensors obtained in Theorems \ref{nw}, \ref{dw} and \ref{Sw}.}
\end{remark}

\begin{proposition}\label{cyclic W} The following identities hold:
\begin{description}
\item [(a)] $\mathfrak{S}_{\beta, \mu, \nu} \ \overcirc{W}^{\alpha}_{\beta\nu\mu} =
\mathfrak{S}_{\beta, \mu, \nu} \ \widehat{W}^{\alpha}_{\beta\nu\mu} =
\mathfrak{S}_{\beta, \mu, \nu} \ R^{a}_{\mu \beta}C^{\alpha}_{\nu a}$
\item [(b)]  $\mathfrak{S}_{\beta, \mu, \nu} \ \widetilde{W}^{\alpha}_{\beta\nu\mu} =
2\mathfrak{S}_{\beta, \mu, \nu} \ R^{a}_{\mu \beta}C^{\alpha}_{\nu a}$
\end{description}
\end{proposition}

\begin{proof} By Theorem \ref{nw}, we have
\begin{eqnarray*} \mathfrak{S}_{\beta, \mu, \nu} \ \overcirc{W}^{\alpha}_{\beta\nu\mu} &=&
\mathfrak{S}_{\beta, \mu, \nu}(\gamma^{\alpha}_{\beta\mu|\nu} - \gamma^{\alpha}_{\beta\nu|\mu})
+ \mathfrak{S}_{\beta, \mu, \nu} (\gamma^{\alpha}_{\beta\epsilon}\Lambda^{\epsilon}_{\nu\mu}) +
\mathfrak{S}_{\beta, \mu, \nu} \ (\gamma^{\epsilon}_{\beta\nu}\gamma^{\alpha}_{\epsilon\mu} -
\gamma^{\epsilon}_{\beta\mu}\gamma^{\alpha}_{\epsilon\nu})\\
&=& \mathfrak{S}_{\beta, \mu, \nu} \ (\Lambda^{\alpha}_{\beta\mu|\nu} + \Lambda^{\epsilon}_{\mu\nu}
\Lambda^{\alpha}_{\beta\epsilon})\\
&=& \mathfrak{S}_{\beta, \mu, \nu} \ (\Lambda^{\alpha}_{\beta\mu|\nu} + \Lambda^{\epsilon}_{\beta\mu}
\Lambda^{\alpha}_{\nu\epsilon} + R^{a}_{\beta\mu}C^{\alpha}_{\nu a}) +
\mathfrak{S}_{\beta, \mu, \nu}R^{a}_{\mu\beta}C^{\alpha}_{\nu a}\\
&=& \mathfrak{S}_{\beta, \mu, \nu}R^{a}_{\mu\beta}C^{\alpha}_{\nu a},\\
\\[- 1cm]\end{eqnarray*}
where in the last step we have used (\ref{FBI}). The proof of the other part of {\bf (a)} is
achieved by applying the first Bianchi identity to the symmetric $d$-connection taking into account Theorem
\ref{Sw} {\bf (a)} together
with the fact that $\widehat{C}^{\alpha}_{\nu a} = C^{\alpha}_{\nu a}$.
The proof of {\bf (b)} is carried out in a similar manner, again by using (\ref{FBI}),
taking into consideration Theorem \ref{dw} {\bf (a)}.
\end{proof}

Summing up, the EAP-space has three distinct $W$-tensors (corresponding to the natural metric, dual and
symmetric $d$-connections), each with six counterparts.
{\bf Eight} only out of the eighteen are {\bf independent},
four coincide with the corresponding curvature tensors and four vanish identically.

\end{mysect}


\begin{mysect}{Cartan-type case}

A {\bf drawback} in the construction of the EAP-space is the fact
that the nonlinear connection is assumed to exist {\it a priori},
independently of the vector fields $\lambda$'s defining the
parallelization. It would be more natural and less arbitrary if the
nonlinear connection were expressed in terms of these vector fields.
In this case, all geometric objects of the EAP-space will be defined
solely in terms of the building blocks of the space. Below, we
impose an extra condition on the canonical $d$-connection, namely,
being of Cartan type. The outcome of such a condition is the
accomplishment of our requirment, besides many other interesting
results.

\bigskip

We first recall the definition of a Cartan type $d$-connection \cite{C}, \cite{GLS}.

\begin{definition} A $d$-connection $D = (\Gamma^{\alpha}_{\mu\nu}, \, \Gamma^{a}_{b\mu}, \, C^{\alpha}_{\mu c}, \, C^{a}_{bc})$ on $TM$
is said to be of {\bf Cartan type}  if
$$D^{h}\!\!\,_{X} {\cal C} = 0; \ \ D^{v}\!\!\,_{X} {\cal C} = v X; \ \ \forall X\in \mathfrak{X}(TM),$$
where ${\cal C} = y^{a}\dot{\partial_a}$ is the Loiuville vector field.
\end{definition}

Locally, the above conditions are expressed in the form
\begin{equation}y^{a}\!\,_{|\mu} = 0, \ \ \ \ y^{a}\!\,_{||c} = \delta^{a}_{c},\end{equation}
or, equivalently,
\begin{equation}N^{a}_{\mu} = y^{b}\Gamma^{a}_{b\mu}, \ \ \ \ y^{b}{C}^{a}_{bc} = 0.\end{equation}

\begin{proposition}\label{TCC} If a $d$-connection $D$ is of Cartan type, then the following identities involving the
torsions and curvatures hold:
\begin{equation}\label{RRTT}R^{a}_{\mu\nu} = y^{b}R^{a}_{b\nu\mu}, \ \ P^{a}_{\mu c} = y^{b}P^{a}_{b\mu c}, \ \
T^{a}_{bc} = y^{d}S^{a}_{dcb}.\end{equation}
\end{proposition}
\begin{proof}
Follows by applying the commutation formulae
{\bf (d), (e)} and {\bf (f)} of Lemma \ref{CFA} to the vector field $y^{a}$.
\end{proof}

\begin{theorem}\label{CT} Let $(TM, \ \lambda)$ be an EAP-space.
Assume that the canonical $d$-connection $D$ is of Cartan type. Then we have:
\begin{description}
\item [(a)] The expression
$y^{b}( \, \undersym{\lambda}{i}^{a}\partial_{\mu} \,
\undersym{\lambda}{i}_{b})$ represents the coefficients of a
nonlinear connection which coincides with the given nonlinear
connection $N^{a}_{\mu}$: $N^{a}_{\mu} =  y^{b}( \,
\undersym{\lambda}{i}^{a}\partial_{\mu} \
\undersym{\lambda}{i}_{b})$.\footnote{A similar expression is found
in \lq\lq ArXiv: 0801.1132 [gr-qc]\rq\rq, but in a completely
different situation.}

\item [(b)] $R^{a}_{\mu\nu} = P^{a}_{\mu c} = T^{a}_{bc} = 0.$ \\Consequently, $C^{a}_{bc}$ is symmetric,
$\gamma^{a}_{bc} = 0$ and $C^{a}_{bc} = \widetilde{C}^{a}_{bc} = \widehat{C}^{a}_{bc} = \ \overcirc{C}^{a}_{bc}.$
\item [(c)] $\lambda^{a}\!\,\,_{\widetilde{||}b} = \lambda^{a}\!\,\,_{\widehat{||}b} =
\lambda^{a}\!\,\,_{{o\atop{||}}b} = 0, \ \ g_{ab\widetilde{||}c} = g_{ab\widehat{||}c} =
g_{ab{o\atop{||}}c} = 0$, \ \ $y^{a}\!\,\,_{\widetilde{||}b} = y^{a}\!\,\,_{\widehat{||}b} =
y^{a}\!\,\,_{{o\atop{||}}b} = \delta^{a}_{b}$.
\item [(d)] $\lambda_{a}$ are positively homogeneous of degree $0$ in $y$. Consequently, so are $g_{ab}$.
\item [(e)] $\dot{\partial_b}N^{a}_{\mu} = \Gamma^{a}_{b\mu}$ and $N^{a}_{\mu}$ is homogeneous.
\\Consequently, $\Gamma^{a}_{b\mu}$ is positively homogeneous of degree $0$ in $y$.
\item [(f)] $\gamma^{a}_{b\mu} = 0$. Consequently, $\Gamma^{a}_{b\mu} = \ \overcirc{\Gamma}^{a}_{b\mu}$.
\item [(g)] $\overcirc{P}^{a}_{\mu b} = 0.$
\end{description}
\end{theorem}

\begin{proof}
We have
\begin{description}
 \item[(a)] $N^{a}_{\mu} = y^{b} (\ \undersym{\lambda}{i}^{a}
\delta_{\mu} \ \undersym{\lambda}{i}_{b}) =  y^{b} \ \undersym{\lambda}{i}^{a}
(\partial_{\mu} - N^{c}_{\mu}\,\dot{\partial_c}) \ \undersym{\lambda}{i}_{b} = y^{b}( \ \undersym{\lambda}{i}^{a}
\partial_{\mu} \ \undersym{\lambda}{i}_{b}) - N^{c}_{\mu} \,y^{b}C^{a}_{bc} =  y^{b}( \ \undersym{\lambda}{i}^{a}
\partial_{\mu} \ \undersym{\lambda}{i}_{b}).$

\item [(b)] is obtained from Proposition \ref{TCC}, taking into account Theorem \ref{v}.
The vanishing of $\gamma^{a}_{bc}$ is
readily obtained by Remark \ref{TC} and $T^{a}_{bc} = 0$.

\item [(c)] is a direct consequence of {\bf (b)}.

\item [(d)] By the symmetry of $C^{a}_{bc}$, we have
$0 = y^{b} C^{a}_{bc} = y^{b}C^{a}_{cb} = \ \undersym{\lambda}{i}^{a}(y^{b}\dot{\partial_b} \ \undersym{\lambda}{i}_{c})$
so that, by (\ref{inverse}), $y^{b}\dot{\partial_b}{\lambda}_{c} = 0$. The result follows from Euler's Theorem.

\item [(e)] follows from the fact that $P^{a}_{\mu b} =
\dot{\partial}_bN^{a}_{\mu} - \Gamma^{a}_{b\mu} = 0$ so that
$y^{b}\dot{\partial}_b N^{a}_{\mu} = y^{b}\Gamma^{a}_{b\mu} = N^{a}_{\mu}$. This could be also deduced from the
expression obtained for $N^{a}_{\mu}$ in {\bf (a)}, taking into account that
$\lambda_{a}$ (hence ${\lambda}^{a}$) are
positively homogeneous of degree $0$ in $y$.

\item [(f)] By {\bf (e)}, we have
$\Gamma^{a}_{b\mu} = \dot{\partial_b}N^{a}_{\mu}$. Consequently,
by definition of the natural metric $d$-connection (Theorem \ref{metric}),
$\frac{1}{2} \ g^{ac}(\delta_{\nu} g_{bc} - g_{dc} \ \dot{\partial_b}N^{d}_{\nu} -
g_{bd} \ \dot{\partial_c}N^{d}_{\nu}) = \ \overcirc{\Gamma}^{a}_{b\mu} - \Gamma^{a}_{b\mu} =
- \ \gamma^{a}_{b\mu}$. Multiplying by $g_{ae}$, we get
$\frac{1}{2} \ (\delta_{\nu} g_{be} - g_{de} \ \dot{\partial_b}N^{d}_{\nu} -
g_{bd} \ \dot{\partial_e}N^{d}_{\nu}) =  - \gamma_{eb\mu}$. This implies that $\gamma_{eb\mu}$ is symmetric in the
indices $e, b$. By Proposition \ref{skew}, $\gamma_{eb\mu}$ is also skew-symmetric in the indices $e, b$.
Consequently, $\gamma_{eb\mu}$ vanishes so that $\gamma^{a}_{b\mu} = g^{ae}\gamma_{eb\mu} = 0$.

\item [(g)] is a direct consequence of the relation
\ $\overcirc{P}^{a}_{\mu b} = P^{a}_{\mu b} + \gamma^{a}_{b\mu}$, taking into account that
$P^{a}_{\mu b} = \gamma^{a}_{b\mu} = 0$.
\end{description}
\vspace{- 0.8 cm}\end{proof}

\begin{corollary}\label{CC} If the canonical $d$-connection is of Cartan type, then
$\widetilde{D}$, $\widehat{D}$ and \ $\overcirc{D}$ are
also of Cartan type.
\end{corollary}

In what follows, we assume that the canonical $d$-connection is of Cartan type.
The next two results are immediate consequences of the fact that
$$R^{a}_{\mu\nu} = P^{a}_{\mu c} = T^{a}_{bc} = \gamma^{a}_{bc} = \gamma^{a}_{b\mu} = 0,$$
taking into consideration Proposition \ref{cyclic W} and Theorems
\ref{n}, \ref{d}, \ref{S}, \ref{nw}, \ref{dw} and \ref{Sw}.

\begin{proposition}\label{WBI} The following relations hold:

\begin{description}
\item [(a)] \ $\overcirc{R}^{a}_{b\mu\nu} = \ \overcirc{P}^{a}_{b\mu c} = \ \overcirc{S}^{a}_{bcd} = 0,$
\item [(b)]$\widetilde{R}^{\alpha}_{\beta\mu\nu} = \Lambda^{\alpha}_{\nu\mu|\beta},$
\item [(c)] $\widetilde{R}^{a}_{b\mu\nu} = \widehat{R}^{a}_{b\mu\nu} =
\widetilde{P}^{a}_{b\mu c} = \widehat{P}^{a}_{b\mu c} = \widetilde{S}^{a}_{bcd} =
\widehat{S}^{a}_{bcd}= 0.$
\item [(d)] \ $\overcirc{W}^{\alpha}_{\beta\nu\mu} = \ \overcirc{R}^{\alpha}_{\beta\mu\nu},$
\item [(e)] \ $\overcirc{W}^{a}_{b\mu\nu} = \ \overcirc{W}^{a}_{b\mu c} = \ \overcirc{W}^{a}_{bcd} =
\widetilde{W}^{a}_{b\mu c} = \widehat{W}^{\alpha}_{b\mu c} = \widetilde{W}^{a}_{bcd} =
\widehat{W}^{a}_{bcd}= 0,$
\item [(f)] $\mathfrak{S}_{\beta, \mu, \nu} \ \overcirc{W}^{\alpha}_{\beta\nu\mu} =
\ \mathfrak{S}_{\beta, \mu, \nu} \ \widehat{W}^{\alpha}_{\beta\nu\mu} =
\ \mathfrak{S}_{\beta, \mu, \nu} \ \widetilde{W}^{\alpha}_{\beta\nu\mu} = 0$.
\end{description}
Consequently, \ $\overcirc{W}^{\alpha}_{\beta\nu\mu}$, $\widehat{W}^{\alpha}_{\beta\nu\mu}$ and
$\widetilde{W}^{\alpha}_{\beta\nu\mu}$ satisfy the first Bianchi identity of the
Riemannian curvature tensor.
\end{proposition}

\begin{theorem}\label{4W} The independent non-vanishing $W$-tensors are given by:
\begin{description}
\item [(a)] $\overcirc{W}^{\alpha}_{\beta\nu c} = (\gamma^{\alpha}_{\beta c|\nu} -
\gamma^{\alpha}_{\beta\nu||c})  + (\gamma^{\epsilon}_{\beta\nu}\gamma^{\alpha}_{\epsilon c} -
\gamma^{\epsilon}_{\beta c}\gamma^{\alpha}_{\epsilon\nu}) -
\gamma^{\epsilon}_{\nu c}\gamma^{\alpha}_{\beta\epsilon}$
\item [(b)] $\widetilde{W}^{\alpha}_{\beta\nu\mu} = \Lambda^{\alpha}_{\nu\mu|\beta} +
\Lambda^{\epsilon}_{\nu\mu}\Lambda^{\alpha}_{\beta\epsilon},$
\item [(c)] $\widetilde{W}^{\alpha}_{\beta\mu c} = \Lambda^{\alpha}_{\mu\beta||c}$
\end{description}
\end{theorem}

To sum up, the assumption that the canonical $d$-connection being of Cartan type implies
that all the geometric objects defined in the EAP-space are expressed in terms
of the vector fields ${\lambda}$'s only. The curvature of the nonlinear connection vanishes and the
three other defined $d$-connections
of the EAP-space, namely the dual, symmetric and the natural metric $d$-connections, are also of Cartan type.
 \ Moreover, there
are only seven non-vanishing curvature tensors and only three independent non-vanishing $W$-tensors,
some of which have simpler expressions than that obtained in the general case.
The EAP-space becomes richer as new relations
among its various geometric objects - which are not valid in the general case - emerge.
Accordingly, the EAP-space in this case becomes more tangible, thus more suitable for physical applications.

\bigskip

We end this section with the following table.

\begin{center} {\bf Table 3: EAP-geometry under the Cartan type case}\\[0.2 cm]
\footnotesize{\begin{tabular}
{|c|c|c|c|c|c|c|c|c|c|c|c|}\hline
&&&\\
{\bf Connection}&{\bf Coefficients}&{\bf Torsion}&{\bf Curvature}\\
 [0.2cm]\hline
&&&\\
{\bf Canonical} &\footnotesize${(\Gamma^{\alpha}_{\mu\nu}, \ \dot{\partial_b}N^{a}_{\nu},
\ C^{\alpha}_{\mu c}, \ C^{a}_{bc})}$&\footnotesize${(\Lambda^{\alpha}_{\mu\, \nu}, \ 0,
\ C^{\alpha}_{\mu c}, \ 0, \ 0)}$&$(0, \ 0, \ 0, \ 0, \ 0, \ 0)$\\[0.2cm]\hline
 &&&\\
{\bf Dual} &\footnotesize${(\Gamma^{\alpha}_{\nu\mu}, \ \dot{\partial_b}N^{a}_{\nu},
\ C^{\alpha}_{\mu c}, \ C^{a}_{bc})}$&\footnotesize${(- \Lambda^{\alpha}_{\mu\nu},
\ 0, \ C^{\alpha}_{\mu c}, \ 0, \ 0)}$&$(\widetilde{R}^{\alpha}_{\beta\mu\nu}, \ 0,
\ \widetilde{P}^{\alpha}_{\beta\mu c}, \ 0, \ 0, \ 0)$\\[0.2cm]\hline
 &&&\\
{\bf Symmetric} &\footnotesize${(\Gamma^{\alpha}_{(\mu\nu)}, \ \dot{\partial_b}N^{a}_{\nu}, \
C^{\alpha}_{\mu c}, \ C^{a}_{bc})}$
&\footnotesize${ (0,
\ 0, \ C^{\alpha}_{\mu c}, \ 0, \ 0)}$&$(\widehat{R}^{\alpha}_{\beta\mu\nu}, \ 0,
\ \widehat{P}^{\alpha}_{\beta\mu c}, \ 0, \ 0, \ 0)$\\[0.2cm]\hline
&&&\\
{\bf Natural} &\footnotesize${( \ \overcirc{\Gamma}^{\alpha}_{\mu\nu},
\ \dot{\partial_b}N^{a}_{\nu}, \ \, \overcirc{C}^{\alpha}_{\mu c}, \ {C}^{a}_{bc})}
$&
\footnotesize${(0, \ 0, \ \, \overcirc{C}^{\alpha}_{\mu c},
\ 0, \ 0)}$&$( \ \overcirc{R}^{\alpha}_{\beta\mu\nu}, \ 0, \ \, \overcirc{P}^{\alpha}_{\beta\mu c}, \ 0,
\ \, \overcirc{S}^{\alpha}_{\beta cd}, \ 0)$
\\[0.2 cm]\hline

\end{tabular}}
\end{center}

\end{mysect}


\begin{mysect}{Berwald-type case}

In this section we assume that the canonical $d$-connection $D$ is of Berwald type
\linebreak \cite{GLS}, \cite{M}. The consequences
of this assumption are investigated.

\begin{definition} A $d$-connection $D = (\Gamma^{\alpha}_{\mu\nu}, \, \Gamma^{a}_{b\mu}, \, C^{\alpha}_{\mu c}, \, C^{a}_{bc})$ on $TM$
is said to be of {\bf Berwald type} if
\begin{equation} \ \ \dot{\partial_{b}}N^{a}_{\mu} = \Gamma^{a}_{b\mu}; \ \ \ \ \ \ \ \ \ C^{\alpha}_{\mu c} = 0.\end{equation}
\\[- 1.5 cm]\end{definition}

\begin{proposition} \label{CCB}If the canonical $d$-connection $D$ is of Cartan type such that $C^{\alpha}_{\mu c} = 0$,
then it is of Berwald type.
\end{proposition}
\begin{proof} Follows from the fact that $0 = P^{a}_{\mu b} = \dot{\partial_b}N^{a}_{\mu} - \Gamma^{a}_{b\mu}$.
\end{proof}

\begin{theorem} \label{BT} Let $(TM, \ \lambda)$ be an EAP-space. Assume that the canonical $d$-connection $D$ is of
Berwald type. Then we have:
\begin{description}
\item [(a)] $P^{a}_{\mu b} = 0$
\item [(b)] ${\lambda}_{\mu}$ are functions of the positional argument $x$ only.
Consequently, so are $g_{\mu\nu}$.
\item [(c)] $ \ \overcirc{C}^{\alpha}_{\mu c} = 0$. Consequently, $\gamma^{\alpha}_{\mu c} = 0$.
\item [(d)] The coefficients $\Gamma^{\alpha}_{\mu\nu}$ and \ $\overcirc{\Gamma}^{\alpha}_{\mu\nu}$ are
functions of the positional argument $x$ only and are given respectively by
$$\Gamma^{\alpha}_{\mu\nu}(x) = ( \ \undersym{\lambda}{i}^{\alpha}\partial_{\nu} \
\undersym{\lambda}{i}_{\mu})(x), \ \ \ \ \ \ \overcirc{\Gamma}^{\alpha}_{\mu\nu}(x) = \frac{1}{2} \, g^{\alpha\epsilon}
(\partial_{\mu}g_{\nu\epsilon} + \partial_{\nu}g_{\mu\epsilon} - \partial_{\epsilon}g_{\mu\nu})(x).$$
\item [(e)] $\Lambda^{\alpha}_{\mu\nu}$ and $\gamma^{\alpha}_{\mu\nu}$ are functions of the positional argument
$x$ only.

\item [(f)] $\gamma^{a}_{b\mu} = \ \overcirc{P}^{a}_{b\mu} = 0.$ Consequently, $\Gamma^{a}_{b\mu} = \
\overcirc{\Gamma}^{a}_{b\mu}.$
\end{description}
\end{theorem}

\begin{proof} The proof is straightforward except for the relation $\gamma^{a}_{b\mu} = 0$,
which can be proved in exactly the same manner as {\bf (f)} of
Theorem \ref {CT}.
\end{proof}

\begin{corollary}\label{BB} If the canonical $d$-connection $D$ is of Berwald type, then
$\widetilde{D}$, $\widehat{D}$ and \ $\overcirc{D}$ are also of Berwald type.
\end{corollary}

In what follows, we assume that the canonical $d$-connection is of Berwald type.
The next two results are immediate consequences of the fact that
$$P^{a}_{\mu c} = C^{\alpha}_{\mu c} = \gamma^{\alpha}_{\mu c} = \gamma^{a}_{b\mu} = 0,$$
taking into account Theorems \ref{n}, \ref{d}, \ref{S}, \ref{nw}, \ref{dw} and \ref{Sw}.

\begin{proposition} The following relations hold:

\begin{description}
\item [(a)] \ $\overcirc{R}^{a}_{b\mu\nu} = \gamma^{a}_{bd}R^{d}_{\mu\nu}, \ \ \ \
\overcirc{P}^{\alpha}_{\beta\nu  c} = \ \overcirc{S}^{\alpha}_{\beta cd} = 0,$
\item [(b)] $\widetilde{R}^{\alpha}_{\beta\nu\mu} = \Lambda^{\alpha}_{\mu\nu|\beta}, \ \ \ \
\widetilde{P}^{\alpha}_{\beta\mu c} = \widehat{P}^{\alpha}_{\beta\mu c} = 0, \ \ \ \
\widetilde{P}^{a}_{b\mu c} = \widetilde{W}^{a}_{b\mu c},$
\item [(c)] $\overcirc{W}^{a}_{b\nu\mu} = \ \overcirc{W}^{\alpha}_{\beta\nu c} =
\widetilde{W}^{\alpha}_{\beta\nu c} = \widehat{W}^{\alpha}_{\beta\nu c} = 0.$

\end{description}
\end{proposition}

\begin{theorem} The independent non-vanishing $W$-tensors are given by:
\begin{description}
\item [(a)] \ $\overcirc{W}^{\alpha}_{\beta\nu\mu} = (\gamma^{\alpha}_{\beta\mu|\nu} -
\gamma^{\alpha}_{\beta\nu|\mu}) +
(\gamma^{\epsilon}_{\beta\nu}\gamma^{\alpha}_{\epsilon\mu} -
\gamma^{\epsilon}_{\beta\mu}\gamma^{\alpha}_{\epsilon\nu}) - \gamma^{\alpha}_{\beta\epsilon}
\Lambda^{\epsilon}_{\nu\mu},$
\item [(b)] \ $\overcirc{W}^{a}_{b\nu c} = \gamma^{a}_{bc|\nu},$
\item [(c)] $\widetilde{W}^{\alpha}_{\beta\nu\mu} = \Lambda^{\alpha}_{\nu\mu|\beta} + \Lambda^{\epsilon}_{\nu\mu}
\Lambda^{\alpha}_{\beta\epsilon},$
\item [(d)] $\widetilde{W}^{a}_{bdc} = T^{a}_{dc||b} + T^{e}_{dc}T^{a}_{be}.$
\end{description}
\end{theorem}

\vspace{0.1 cm}To sum up, the assumption that the canonical $d$-connection being of Berwald type implies that most of the purely
horizontal geometric objects of the EAP-space become functions of the positional argument $x$ only and
coincide with the corresponding geometric objects of the classical AP-space. Moreover, the three other
defined $d$-connections turn out also to be of Berwald type. Finally, in this case, there are twelve non-vanishing
curvature tensors and four independent $W$-tensors.

\bigskip

We end this section with the following table (compare with Table 3).

\begin{center} {\bf Table 4: EAP-geometry under the Berwald type case}\\[0.2cm]
\footnotesize{\begin{tabular}
{|c|c|c|c|c|c|c|c|c|c|c|c|}\hline
&&&\\
{\bf Connection}&{\bf Coefficients}&{\bf Torsion}&{\bf Curvature}\\
 [0.2cm]\hline
&&&\\
{\bf Canonical} &\footnotesize${(\Gamma^{\alpha}_{\mu\nu}, \ \dot{\partial_b}N^{a}_{\nu},
\ 0, \ C^{a}_{bc})}$&\footnotesize${(\Lambda^{\alpha}_{\mu\, \nu}, \ R^{a}_{\mu\nu},
\ 0, \ 0, \ T^{a}_{bc})}$&$(0, \ 0, \ 0, \ 0, \ 0, \ 0)$\\[0.2cm]\hline
 &&&\\
{\bf Dual} &\footnotesize${(\Gamma^{\alpha}_{\nu\mu}, \ \dot{\partial_b}N^{a}_{\nu},
\ 0, \ C^{a}_{cb})}$&\footnotesize${(- \Lambda^{\alpha}_{\mu\nu},
\ R^{a}_{\mu\nu}, \ 0, \ 0, \ - T^{a}_{bc})}$&$(\widetilde{R}^{\alpha}_{\beta\mu\nu},
\,\widetilde{R}^{a}_{b\mu\nu}, \, 0,
\ \widetilde{P}^{\alpha}_{\beta\mu c}, \, 0, \, \widetilde{S}^{a}_{bcd})$\\[0.2cm]\hline
 &&&\\
{\bf Symmetric} &\footnotesize${(\Gamma^{\alpha}_{(\mu\nu)}, \ \dot{\partial_b}N^{a}_{\nu}, \
0, \ C^{a}_{(bc)})}$
&\footnotesize${(0,
\ R^{a}_{\mu\nu}, \ 0, \ 0, \ 0)}$&$(\widehat{R}^{\alpha}_{\beta\mu\nu}, \, \widehat{R}^{a}_{b\mu\nu},
\, 0, \, \widehat{P}^{\alpha}_{\beta\mu c}, \, 0, \, \widehat{S}^{a}_{bcd})$\\[0.2cm]\hline
&&&\\
{\bf Natural} &\footnotesize${( \ \overcirc{\Gamma}^{\alpha}_{\mu\nu},
\ \dot{\partial_b}N^{a}_{\nu}, \ 0, \ \, \overcirc{C}^{a}_{bc})}
$&
\footnotesize${(0, \ R^{a}_{\mu\nu}, \ 0,
\ 0, \ 0)}$&$( \ \overcirc{R}^{\alpha}_{\beta\mu\nu}, \ \overcirc{R}^{a}_{b\mu\nu}, \, 0, \
\overcirc{P}^{a}_{b\mu c}, \, 0, \ \overcirc{S}^{a}_{bcd})$
\\[0.2 cm]\hline

\end{tabular}}
\end{center}

\end{mysect}

\begin{mysect}{Recovering the classical AP-space}

We now assume that the canonical $d$-connection $D$ is both Cartan and Berwald type. In view of
Proposition \ref{CCB}, this condition is equivalent to the (apparently weaker) condition
that $D$ is of Cartan type and $C^{\alpha}_{\mu c} = 0$. We show that
in this case the classical AP-space emerges, in a natural way, as a special case from the EAP-space.
We refer to this condition as the {\bf CB}-condition.

$$\text{{\bf CB}-condition:} \ N^{a}_{\mu} = y^{b} \Gamma^{a}_{b\mu}, \ \ \ \
y^{b}C^{a}_{bc} = 0; \ \ \ \ C^{\alpha}_{\mu c} = 0.$$

As easily checked, we have

\begin{theorem}\label{CartanB} Assume that the {\bf CB}-condition holds. Then
\begin{description}
\item [(1)] The four defined $d$-connections of the EAP-space coincide up to the $hh$-coefficients. Moreover, these
$hh$-coefficients are functions of the
positional argument $x$ only and are identical to the coefficients of the corresponding
four defined connections in the
classical AP-space.
\item [(2)] The torsion of the canonical $d$-connection and the contortion of the EAP-space are functions
of the positional argument $x$ only and are given by
$${\bf T} = (\Lambda^{\alpha}_{\mu\nu}, \ 0, \ 0, \ 0 , \ 0); \ \ \ \ {\bf C} =
(\gamma^{\alpha}_{\mu\nu}, \ 0, \ 0, \ 0)$$
\item [(3)] The three non-vanishing curvature tensors are functions of the positional \linebreak argument $x$ only and
are given by
\begin{description}
\item [(a)] \!\ $\overcirc{R}^{\alpha}_{\beta\mu\nu} = (\gamma^{\alpha}_{\beta\mu|\nu} -
\gamma^{\alpha}_{\beta\nu|\mu}) + (\gamma^{\epsilon}_{\beta\nu}\gamma^{\alpha}_{\epsilon\mu} -
\gamma^{\epsilon}_{\beta\mu}\gamma^{\alpha}_{\epsilon\nu}) + \gamma^{\alpha}_{\beta\epsilon}
\Lambda^{\epsilon}_{\mu\nu}$
\item [(b)] $\widetilde{R}^{\alpha}_{\beta\mu\nu} = \Lambda^{\alpha}_{\nu\mu|\beta},$
\item [(c)] $\widehat{R}^{\alpha}_{\beta\mu\nu} =
\frac{1}{2}\, (\Lambda^{\alpha}_{\beta\mu|\nu} - \Lambda^{\alpha}_{\beta\nu|\mu}) +
\frac{1}{4}\, (\Lambda^{\epsilon}_{\beta\mu} \Lambda^{\alpha}_{\nu\epsilon} -
\Lambda^{\epsilon}_{\beta\nu} \Lambda^{\alpha}_{\mu\epsilon}) +
\frac{1}{2}\, (\Lambda^{\epsilon}_{\mu\nu}\Lambda^{\alpha}_{\beta\epsilon})$
\end{description}
\item [(4)] There is only one $W$-tensor which is a function of the positional
argument $x$ only and is given
by
$$\widetilde{W}^{\alpha}_{\beta\nu\mu} = \Lambda^{\alpha}_{\nu\mu|\beta} + \Lambda^{\epsilon}_{\nu\mu}
\Lambda^{\alpha}_{\beta\epsilon}$$
All other W-tensors vanish identically, or coincide with the corresponding curvature tensors.
\end{description}
Consequently, the fundamental geometric objects of the EAP-space coincide with the corresponding
fundamental geometric objects of the classical AP-space \cite{AMR}.
\end{theorem}

\begin{corollary}\label{appc} If the canonical $d$-connection $D$ satisfies the {\bf CB}-condition, then
$\widetilde{D}$, $\widehat{D}$ and \ $\overcirc{D}$ also satisfy the {\bf CB}-condition.
\end{corollary}

We end this section with the following two tables which summarize the geometry of the EAP-space
under the {\bf CB}-condition.

\begin{center} {\bf Table 5: Fundamental connections of EAP-space \\under the CB-condition}\\[0.3cm]
\footnotesize{\begin{tabular}
{|c|c|c|c|c|c|c|c|c|c|c|c|}\hline
&&\\
{\bf Connection}&{\bf Coefficients}&{\bf hh-Coefficients}\\[0.2cm]\hline
&&\\
{\bf Canonical}&\footnotesize${(\Gamma^{\alpha}_{\mu\nu}, \ \dot{\partial_{b}}N^{a}_{\nu},
\ 0, \ C^{a}_{bc})}$& $\Gamma^{\alpha}_{\mu\nu}(x) =
( \ \undersym{\lambda}{i}^{\alpha}\partial_{\nu} \ \undersym{\lambda}{i}_{\mu})(x)$\\[0.2cm]\hline
&&\\
{\bf Dual}&\footnotesize${(\widetilde{\Gamma}^{\alpha}_{\mu\nu}, \ \dot{\partial_{b}}N^{a}_{\nu},
\ 0, \ C^{a}_{bc})}$& $\widetilde{\Gamma}^{\alpha}_{\mu\nu}(x) =
{\Gamma}^{\alpha}_{\nu\mu}(x)$\\[0.2cm]\hline
&&\\
{\bf Symmetric}&\footnotesize${(\widehat{\Gamma}^{\alpha}_{\mu\nu}, \ \dot{\partial_{b}}N^{a}_{\nu}, \ 0, \ C^{a}_{bc})}$
& $\widehat{\Gamma}^{\alpha}_{\mu\nu}(x) = \Gamma^{\alpha}_{(\mu\nu)}(x)$\\[0.2cm]\hline
&&\\
{\bf Natural}&\footnotesize${( \ \overcirc{\Gamma}^{\alpha}_{\mu\nu},
\ \dot{\partial_{b}}N^{a}_{\nu},  \ 0, \ C^{a}_{bc})}
$& \ $\overcirc{\Gamma}^{\alpha}_{\mu\nu}(x) = \frac{1}{2} \
g^{\alpha\epsilon}(\partial_{\mu}g_{\nu\epsilon} +
\partial_{\nu}g_{\mu\epsilon} - \partial_{\epsilon}g_{\mu\nu})(x)$\\[0.2 cm]\hline
\end{tabular}}
\end{center}

\vspace{0.05 cm}

\begin{center} {\bf Table 6: EAP-geometry under the CB-condition}\\[0.2 cm]
\footnotesize{\begin{tabular}
{|c|c|c|c|c|c|c|c|c|c|c|c|}\hline
&&&\\
{\bf Connection}&{\bf Coefficients}&{\bf Torsion}&{\bf Curvature}\\
 [0.2cm]\hline
&&&\\
{\bf Canonical} &\footnotesize${(\Gamma^{\alpha}_{\mu\nu}, \ \dot{\partial_b}N^{a}_{\nu},
\ 0, \ C^{a}_{bc})}$&\footnotesize${(\Lambda^{\alpha}_{\mu\, \nu}, \ 0,
\ 0, \ 0, \ 0)}$&$(0, \ 0, \ 0, \ 0, \ 0, \ 0)$\\[0.2cm]\hline
 &&&\\
{\bf Dual} &\footnotesize${(\Gamma^{\alpha}_{\nu\mu}, \ \dot{\partial_b}N^{a}_{\nu},
\ 0, \ C^{a}_{bc})}$&\footnotesize${( - \Lambda^{\alpha}_{\mu\nu},
\ 0, \ 0, \ 0, \ 0)}$&$(\widetilde{R}^{\alpha}_{\beta\mu\nu}, \ 0, \ 0, \ 0, \ 0, \ 0)$\\[0.2cm]\hline
 &&&\\
{\bf Symmetric} &\footnotesize${(\Gamma^{\alpha}_{(\mu\nu)}, \ \dot{\partial_b}N^{a}_{\nu}, \
0, \ C^{a}_{bc})}$
&\footnotesize${(0,
\ 0, \ 0, \ 0, \ 0)}$&$(\widehat{R}^{\alpha}_{\beta\mu\nu}, \ 0, \ 0, \ 0, \ 0, \ 0)
$\\[0.2cm]\hline
&&&\\
{\bf Natural} &\footnotesize${( \ \overcirc{\Gamma}^{\alpha}_{\mu\nu},
\ \dot{\partial_b}N^{a}_{\nu}, \ 0, \ {C}^{a}_{bc})}
$&
\footnotesize${(0, \ 0, \ 0,
\ 0, \ 0)}$&$( \ \overcirc{R}^{\alpha}_{\beta\mu\nu}, \ 0, \ 0, \ 0, \ 0, \ 0)$
\\[0.2 cm]\hline

\end{tabular}}
\end{center}

\bigskip
{\bf Some comments on the {\bf CB}-condition:}
\begin{description}

\item [(a)] It should be noted that the non-vanishing of the purely vertical
tensors ${\lambda}^{a}$, $C^{a}_{bc}$, and $g_{ab}$ and the
$hv$-coefficients $\Gamma^{a}_{b\mu}$ of the canonical
$d$-connection may represent extra degrees of freedom which do not
exist in the classical AP-context. Moreover, these vertical
geometric objects still depend on the directional argument $y$.
However, they actually {\it do not contribute to the EAP-geometry
under the {\bf CB}-condition}. This is because the torsion and the
contortion tensors in this case have only one non-vanishing
counterpart each, namely the purely horizontal components
$\Lambda^{\alpha}_{\mu\nu}$ and $\gamma^{\alpha}_{\mu\nu}$
respectively.

\item [(b)] One reading of Theorem \ref{CartanB} is roughly that the \lq\lq projection\rq\rq\, of the geometric objects of the
EAP-space on the base manifold
$M$ can be identified with the classical AP-geometry.
The distinction
which appears between the two geometries is due to the fact that the geometric objects of the EAP-space
live in the double tangent bundle
$TTM\to TM$ and not in the
tangent bundle $TM\to M$.
Consequently, it \linebreak can be said, roughly speaking, {\it that the classical AP-space is a copy of the
EAP-space equipped with the {\bf CB}-condition, viewed from the perspective of the \linebreak base manifold $M$}.
\end{description}
\end{mysect}

\begin{mysect}{Concluding remarks}

In the present article, we have constructed and developed a parallelizable structure analogous to the AP-geometry on the tangent bundle $TM$ of $M$.
Four linear connections, depending on one {\it a priori} given nonlinear connection, are used to explore the properties of this geometry.
Different curvature
tensors charaterizing this structure, together with their contractions, are computed. The different $W$-tensors are also derived. Extra conditions
are imposed on the canonical $d$-connection, the consequences of which are investigated. Finally, a special form of the canonical $d$-connection is
introduced under which the EAP-geometry
reduces to the classical AP-geometry.

\bigskip

On the present work, we have the following comments:
\begin{description}

\item [(1)] Existing theories of gravity suffer from some problems connected to
recent observed astrophysical phenomena, especially those admitting
{\bf anisotropic} behavior of the system concerned (e.g. the
flatness of the rotation curves of spiral galaxies). So, theories in
which the gravitational potential depends {\it on both position and
direction} may be needed. Such theories are to be constructed in
spaces admitting this dependence; a potential candidate is the
EAP-space. This is one of the aims motivating the present work.

\item [(2)] One possible physical application of the EAP-geometry would be the construction of a generalized field theory
on the tangent bundle $TM$ of $M$. This
could be achieved by a double contraction of the purely horizontal $W$-tensor $W^{\alpha\beta}\,_{\mu\nu}$ and the
purely vertical $W$-tensor $W^{ab}\,_{cd}$
to obtain respectively the \lq\lq horizontal\rq\rq \, scalar Lagrangian
${\cal H} := |\lambda| H := |\lambda| g^{\alpha\beta} H_{\alpha\beta}$, where $|\lambda|: = det(\lambda_{\alpha})$ and
$$H_{\alpha\beta} := \Lambda^{\nu}_{\epsilon\alpha}\Lambda^{\epsilon}_{\nu\beta} - C_{\alpha}C_{\beta}$$
and the \lq\lq vertical\rq\rq \, scalar Lagrangian ${\cal V} := {||\lambda||} V := {||\lambda||} g^{ab} V_{ab}$, where
${||\lambda||}: = det(\lambda_{a})$ and
$$V_{ab} := T^{d}_{ea}T^{e}_{db} - C_{a}C_{b}.$$
The field equations are obtained by
the use of the Euler-Lagrange equations
$$ \ \ \, \ \frac{\partial {\cal H}}{\partial \lambda_{\beta}} - \frac{\partial}{\partial x^{\epsilon}}(\frac{\partial {\cal H}}
{\partial \lambda_{\beta, \epsilon}}) -  \frac{\partial}{\partial y^{e}}(\frac{\partial {\cal H}}
{\partial \lambda_{\beta; \, e}}) = 0, \ \ \ \ \ \text{(horizontal form)}$$
$$\, \!\, \frac{\partial {\cal V}}{\partial \lambda_{b}} - \frac{\partial}{\partial x^{\epsilon}}(\frac{\partial {\cal V}}
{\partial \lambda_{b, \epsilon}}) -  \frac{\partial}{\partial y^{e}}(\frac{\partial {\cal V}}
{\partial \lambda_{b; \, e}}) \,= \,0. \ \ \ \ \ \ \, \text{(vertical form)}$$
The resulting field equations could be compared with those derived by M. Wanas \cite{aa} and R. Miron \cite{Miron}.

\item [(3)] Among the advantages of the classical AP-geometry are the ease in calculations and the
diverse and its thorough applications \cite{b}. For these reasons, we have kept, in this work, as
close as possible to the classical AP-case. However, {\it the extra
degrees of freedom in our EAP-geometry have created an abundance of
geometric objects which have no counterpart in the classical
AP-geometry}. Since the physical meaning of most of the geometric
objects of the classical AP-structure is clear, we hope to attribute
physical meaning to the new geometric objects appearing in the
present work. The physical interpretation of the geometric objects existing in the EAP-space and not in the AP-geometry
may need deeper investigation.  The study of the Cartan type case, due to its simplicity,
may be our first step in tackling the general case.

\item [(4)] In conclusion, we hope that physicits working in AP-geometry would divert their attention to the
study of EAP-geometry and its consequences, due to its wealth,
relative simplicity and its close resemblance (at least in form) to
the classical AP-geometry. We believe that the extra degrees of
freedom offered by the EAP-geometry may give us more insight into
the infrastructure of physical phenomena studied in the context of
classical AP-geometry and thus help us better understand the theory
of general relativity and its connection to other physical theories.

\end{description}
\end{mysect}


\bibliographystyle{plain}

\end{document}